\documentclass[a4paper,12pt]{article}
\usepackage{amsthm,amsfonts, indentfirst, latexsym, bm,graphicx,amsmath,amssymb,mathrsfs,verbatim,hyperref,pdfpages}
\usepackage{caption}
\usepackage{subcaption}
\usepackage{color}
\usepackage[T1]{fontenc}
\usepackage{dsfont}
\usepackage{latexsym}
\usepackage{url}
\usepackage{array}
\usepackage{tikz}
\usepackage{pgfplots}
\usepackage{graphicx}
\usepackage[titletoc,title]{appendix}
\usepackage{bbm} 
\usepackage{algorithm,algorithmic}
\numberwithin{equation}{section}

\usepackage{chngcntr}
\counterwithout{equation}{section} 
\usetikzlibrary{arrows, automata,positioning,calc,shapes,decorations.pathreplacing}
\colorlet{greencolor}{green!50!black}
\colorlet{textcolor}{red}
\colorlet{tancolor}{orange!80!black}
\colorlet{bluecolor}{blue}
\definecolor{mycolor2}{rgb}{0.07843,0.16863,0.54902}

\title{Parameter and State 
	Estimation in Queues
	and Related Stochastic Models: A Bibliography}

\author{ Azam Asanjarani        \and Yoni Nazarathy  
}

\author{Azam Asanjarani\footnote{The University of Auckland.}, 
	Yoni Nazarathy\footnote{The University of Queensland.}
	.}

\date{ }
\begin{document}
	\maketitle
	\begin{abstract}
		
		\noindent
		This is an annotated bibliography on estimation and inference results for queues and related stochastic models. The purpose of this document is to collect and  categorise works in the field, allowing for researchers and practitioners to explore the various types of results that exist.\\
		
		\noindent
		{\bf Our focus is on papers that deal with mathematical queueing models as well as related stochastic models motivated by queues}. \\
		
		
		\noindent
		We attempted to make this  bibliography exhaustive, yet there are possibly some papers that we have missed. As it is updated continuously, additions and comments are welcomed. \\
		
		\noindent
		Note that this bibliography is also a companion to our survey of parameter and state estimation in queues \cite{asanjarani2021survey}.
	\end{abstract}
	
	\vspace{40pt}
	{\scriptsize
		Additional works not mentioned in this bibliography include the following related categories: (i) Methods for parameter estimation of point processes (not involving queues). (ii) Methods for parameter estimation of stochastic matrix analytic models (not involving queues). (iii) Optimal control of queues including bandit problems where state or parameter estimation is not directly considered. (iv) Inference, estimation and tomography of communication networks not directly modelled as queueing networks. (v) Analysis of queues with parameter uncertainty, and robustness not directly involving inference. (vi) Road traffic modelling not involving an explicit congestion or queueing model.  }
	
	\newpage
	
	\newpage
	{
		\parindent0mm
		\section{Chronological order with brief descriptions}
		
		{\bf {\color{mycolor2}1955}}
		\medskip

		{\bf Cox \cite{cox1955statistical}: } An overview paper of queueing theory outlining the philosophy of estimating parameters of input processes vs. performance processes.\\
		
		\medskip
		{\bf  {\color{mycolor2}1957}}
		\medskip
		
		{\bf Benes \cite{benes1957sufficient}:} Transient M/M/$\infty$ full observation over a fixed interval.
		
		{\bf Clarke \cite{clarke1957maximum}:} M/M/1 MLE with full observation. The first paper. Sampling until ``busy time'' reaches a pre-assigned value yields closed-form MLEs.\\

		\medskip
		{\bf  {\color{mycolor2}1961}}
		
		\medskip
		
		{\bf Billingsley \cite{billingsley1961statisticalPaper}:} Book on inference of Markov chains.\\
		
		\medskip
		{\bf  {\color{mycolor2}1965}}
		
		\medskip
		
		{\bf Cox \cite{cox1965some}:}
		An overview paper on parameter estimation,  separate analysis of input and service mechanism and problems connected with the sampling of queueing process.
		
		{\bf Kovalenko \cite{1965Kovalenko}:}  
		On recovering the characteristics of a system from observations of the outgoing flow (in Russian)
		
		{\bf Newell \cite{newell1965approximation}:} A review of approximation methods for queues with application to the fixed-cycle traffic light.
		
		{\bf Wolff \cite{wolff1965problems}:} Large sample theory for birth-death queues.\\

		\medskip
		{\bf  {\color{mycolor2}1966}}
		
		\medskip
		
		{\bf Lilliefors \cite{lilliefors1966some}:} Confidence intervals for standard performance measurements based on parameter error.\\
		
		\medskip
		{\bf  {\color{mycolor2}1967}}
		
		\medskip
		
		{\bf Greenberg \cite{greenberg1967behavior}:} Different ways of determining for how long to observe a stationary M/M/1 queue (e.g. fixed number of arrivals, fixed total observation time, etc.\\

		\medskip
		{\bf  {\color{mycolor2}1968}}
		
		\medskip
		
		{\bf Daley \cite{daley1968monte}:}
		The serial correlation coefficients of queue sizes in a stationary GI/M/1
		queue are studied. 
		
		{\bf Daley \cite{daley1968serial}:} 
		The serial correlation coefficients
		of a (stationary) sequence of waiting times in a stationary M/M/1, M/G/1 and GI/G/l queueing system are studied.\\

		\medskip
		{\bf  {\color{mycolor2}1970}}
		
		\medskip
		
		{\bf Brown \cite{brown1970m}:} Estimating the $G$ in M/G/$\infty$ with arrival and departure times without knowing what customers they related to.
		
		{\bf Ross \cite{ross1970identifiability}:} Discusses identifying the distributions of GI/G/k uniquely based on observation of the queueing process.\\
		
		\medskip
		{\bf  {\color{mycolor2}1971}}
		
		\medskip
		{\bf Pakes \cite{pakes1971serial}:} 
		The serial correlation coefficients of waiting times in the stationary GI/M/1 queue are studied. (completing \cite{daley1968monte} work)\\
		
		\medskip
		{\bf  {\color{mycolor2}1972}}
		
		\medskip
		
		{\bf Goyal and Harris \cite{goyal1972maximum}:} MLE for queues with Poisson arrivals with state-dependent general service times when queue sizes are observed at departure points.
		
		{\bf Jenkins \cite{jenkins1972relative}:}  Compares the asymptotic variance of two estimators for M/M/1.

		{\bf Muddapur \cite{muddapur1972bayesian}:} Adds a prior distribution to Clarke's \cite{clarke1957maximum} approach.\\

		\medskip
		{\bf  {\color{mycolor2}1973}}
		
		\medskip
		{\bf Harris  \cite{harris1973some}:} 
		An overview paper presented the statistical analysis of queueing systems with an emphasis on the estimation of input and service parameters and/or distributions. 
		
		{\bf Neal and  Kuczura \cite{neal1973theory}:}
		Presents accurate approximation and asymptotic approximations (by using renewal theory) for the variance of any differentiable functions of different traffic measurements. 
		
		{\bf Reynolds \cite{reynolds1973estimating}:} Bayesian approach for estimation of birth death parameters.\\

		\medskip
		{\bf  {\color{mycolor2}1974}}
		
		\medskip
		
		{\bf Aigner \cite{aigner1974parameter}:} Compares properties of various estimators for M/M/1 with cross-sectional data.
		
		{\bf Brillinger \cite{brillinger1974cross}:} Estimates parameter for a linear time invariant model that generalizes the G/G/$\infty$ queue.\\
		
		\medskip
		{\bf  {\color{mycolor2}1975}}
		
		\medskip
		
		{\bf Keiding \cite{keiding1975maximum}:} Analyses asymptotic properties of the MLE for a birth-and-death process with linear rates (both birth and death).\\

		\medskip
		{\bf  {\color{mycolor2}1979}}
		
		\medskip
		
		{\bf Thiagarajan and Harris \cite{thiagarajan1979statistical}:}  Exponential goodness of fit test for the service times of M/G/1 based on observations of the queue lengths and/or waiting times.\\

		\medskip
		{\bf  {\color{mycolor2}1980}}
		
		\medskip
		
		{\bf Dave and Shah \cite{dave1980maximum}:} MLE of an M/M/2 queue with heterogenous servers.

		{\bf Gordon and Dowdy \cite{gordon1980impact}:} Investigates the effect of parameter estimation errors on the performance in closed product form queueing networks. \\
		
	\newpage
		\medskip
		{\bf  {\color{mycolor2}1981}}
		
		\medskip
		
		{\bf Basawa and Prabhu \cite{basawa1981estimation}:} Estimates for non-parametric and parametric models of single server queues over a horizon up to the nth departure epoch. Also the m.l.e's of the mean inter-arrival time and mean service time in an M/M/1 observed over a fixed time-interval.
		
		{\bf Walrand \cite{walrand1981filtering}:} Proposes an elementary justification of the filtering formulas for a Markov chain and an analysis of the arrival and departure processes at a ./M/1 queue in a quasi-reversible network. 
		
		{\bf Grassmann \cite{grassmann1981technical}:}
		This paper shows that in the M/D/$\infty$ queueing system with service time $S$,
		the optimal way to estimate the expected number in the system is by sampling
		the system at time $0, S, 2S, \cdots , kS$.\\
		
		\medskip
		{\bf  {\color{mycolor2}1982}}
		
		\medskip
		
		{\bf Halfin\cite{halfin1982linear}:} 
		Finds the minimum-variance linear estimator for the expected value of a
		stationary stochastic process observed over a finite time interval, whose
		covariance function is a sum of decaying exponentials.

		{\bf Schruben and Kulkarni \cite{schruben1982some}:} Studies the interface between stochastic models and actual systems for the special case of  M/M/1 queue.\\

		\medskip
		{\bf  {\color{mycolor2}1983}}
		
		\medskip
		
		{\bf Hern{\`a}ndez-Lerma and Marcus \cite{hernandez1983adaptive}:} Adaptive control of an M/G/1 queueing system, where the control chooses the service rate to minimize costs.\\
		
		\medskip
		{\bf  {\color{mycolor2}1984}}
		
		\medskip
		
		{\bf Baras, Dorsey, and Makowski \cite{baras1984stability}:} Stability, parameter estimation and adaptive control for discrete-time competing queues.
		
		{\bf Eschenbach \cite{eschenbach1984statistical}:} 
		Briefly describes methods and results in the statistical analysis of queueing systems.

		{\bf Edelman and McKellar \cite{edelman1984comments}:} Comments on \cite{dave1980maximum}.
		
		{\bf Hern{\`a}ndez-Lerma and Marcus \cite{hernandez1984optimal}:} Adaptive control of priority assignment in a multi-class queue.
		
		{\bf Machihara \cite{machihara1984carried}:}  The carried traffic estimate errors for delay system models are analyzed with emphasis on the analysis of the effect of the holding time distribution on the estimate errors.
		
		{\bf Warfield and Foers \cite{warfield1984application}:} A Bayesian method for analysing teletraffic
		measurement data is discussed.
		
		{\bf Woodside, Stanford and Pagurek \cite{woodside1984optimal}:}
		Presents optimal mean square predictors for queue lengths and delays in
		the stationary GI/M/m queue, based on a queue length measurement.\\
		
		\medskip
		{\bf  {\color{mycolor2}1985}}
		
		\medskip
		
		{\bf Armero \cite{armero1985bayesian}:} 
		The posterior distribution of traffic intensity and the posterior predictive distribution of the waiting time and number of customers for a M/M/1/$\infty$ FIFO queue are obtained given two independent samples of arrival and service times.
		
		{\bf Warfield and Foers \cite{warfield1985application}:} Bayesian analysis for traffic intensity in M/M/c/K type models and in retrial models.\\

		\medskip
		{\bf  {\color{mycolor2}1986}}
		
		\medskip
		
		{\bf Subba Rao and Harishchandra \cite{subba1986large}:} Large normal approximation test based for the traffic intensity parameter in GI/G/s queues.\\
		
		\medskip
		{\bf  {\color{mycolor2}1987}}
		
		\medskip
		{\bf Bhat and Rao \cite{bhat1987statistical}:} A first major survey on statistical analysis of queueing systems. 
		
		\textbf {Mcgrath, Gross and Singpurwalla \cite{mcgrath1987subjective}:}  Attempts to illustrate Bayesian approach through M/M/1 and M/M/1/K examples.
		
		\textbf{ McGrath and Singpurwalla \cite{mcgrath1987subjectiveII}:} This is part II to \cite{mcgrath1987subjective} (without Gross). Here the focus is on integrating the ``Shannon measure of information''(cross-entropy) in the analysis.

		{\bf Ramalhoto \cite{ramalhoto1987some}:} Discusses estimation of generalizations of GI/G/$\infty$, i.e. random translations whose distribution is parameterized by a certain function, $h(\cdot)$.\\
		
		\medskip
		{\bf  {\color{mycolor2}1988}}
		
		\medskip
		
		{\bf Basawa and Prabhu \cite{basawa1988large}:} Estimation of GI/G/1 with exponential family densities. Full observation over $[0,T]$ where $T$ is a stopping time. Several $T$'s considered and asymptotic properties compared.
		
		{\bf Chen, Harrison, Mandelbaum, Van Ackere,  Wein \cite{chen1988empirical}}: Empirical evaluation of a queueing network model for semiconductor wafer fabrication.

		{\bf Harishchandra and Rao \cite{harishchandra1988note}:} Inference for the M/$E_k$/1 queue.

		{\bf Jain and Templeton \cite{jain1988statistical}:} Estimation of GI/M/1 (and GI/M/1/m with m known) parameters where the arrival rate is either $\lambda$ or $\lambda_1$ depending on the queue level.

		{\bf Nozari and Whitt \cite{nozari1988estimating}:} Propose an indirect approach to estimate average production intervals (the length of time between starting and finishing work on each product) using work-in-process inventory measurements. \\
		
		\medskip
		{\bf  {\color{mycolor2}1989}}
		
		\medskip
		
		{\bf Fendick and Whitt \cite{fendick1989measurements}:}Proposes measurements and approximations to describe the variability of offered traffic to a queue (the variability of the arrival process together with the service requirements) and predicts the average workload in the queue (which assumed to have a single server, unlimited waiting space and a work-conserving service discipline).

		{\bf Glynn and Whitt \cite{glynn1989indirect}:} Using the little's $L=\lambda W$ and generalizations to infer $L$ from $W$ and the oppositive.

		{\bf Hantler and Rosberg \cite{hantler1989optimal}:} Parameter estimation of M/M/c queue with parameters in stochastic varying environment, first doing the constant invariant derivation and then using in conjunction with Kalman filter for the time-varying case.
		
		{\bf Jain and Templeton \cite{jain1989problem}:} Sequential analysis view for M/$E_k$/1 queues.
		
		\textbf{Sengupta \cite{sengupta1989markov}:} Present an algorithm for computing the steady-state distribution of the waiting time and queue length of the stable GI/K/l queue.\\
		
		\medskip
		{\bf  {\color{mycolor2}1990}}
		
		\medskip
		
		{\bf Gaver and Jacobs \cite{gaver1990inference}:} Transient M/G/1 inference.
		
		{\bf Gawlick \cite{gawlick1990estimating}:} Applies the (QIE) to ethernet data.
		
		{\bf Larson \cite{larson1990queue}:}  This deals with ``State Reconstruction'' as opposed to ``parameter inference'' in what is called the ``Queue Inference Engine'' (QIE). This is the first of many papers on the idea of using transactional data to reconstruct an estimate of the queue length process.

		{\bf Rubin and Robson \cite{rubin1990single}:} Inference and estimation of number of arrivals for a queueing system with losses due to bulking and a server that works a fixed shift and stays to work after the shift. Small sample analysis as opposed to asymptotic properties.\\
		
		\medskip
		{\bf  {\color{mycolor2}1991}}
		
		\medskip

		{\bf Hall and Larson \cite{hall1991using}:} Modifies (extends) the QIE \cite{larson1990queue} to finite queues and to a case where there is data about exceeding a certain level.
		
		{\bf Jain and Templeton \cite{jain1991confidence}:}  Confidence intervals for estimation for M/M/2 with heterogenous servers.
		
		{\bf  Jain \cite{jain1991comparison}:} Compares confidence intervals for $\rho$ using several methods and sampling regimens in M/$E_k$/1 queues.
		
		{\bf Larson \cite{larson1991queue}:} An addendum to \cite{larson1990queue} reducing the computational complexity from $O(N^5)$ to $O(N^3)$.
		
		{\bf Thiruvaiyaru, Basawa  and Bhat \cite{thiruvaiyaru1991estimation}:} Large sample theory for MLEs of Jackson networks.\\

		\medskip
		{\bf  {\color{mycolor2}1992}}
		
		\medskip
		
		\textbf{Asmussen \cite{asmussen1992phase}:} proves that the stationary waiting time in a GI/PH/1 queue with phase-type service time is phase-type.
		
		\textbf{Asmussen and Bladt \cite{asmussen1992renewal}:} The Matrix-exponential distribution is introduced and some of its basic structural properties are given. Further, an algorithm for computing the waiting time distribution of a queue with matrix-exponential service times and general inter-arrival times is given. This algorithm is a slight generalization of the algorithm for computing the waiting time distribution of $GI/PH/1$ queues. 
		
		{\bf Basawa and Bhat \cite{basawa1992sequential}:} Presents sequential analysis methods for the traffic intensity of single server queues.
		
		{\bf Bertsimas and Servi \cite{bertsimas1992deducing}:} Improves on the $O(n^5)$ algorithm in \cite{larson1990queue} to $O(n^3)$. Also presents an on-line algorithm for estimating the queue length after each departure and includes time-varying Poisson generalizations.
		
		{\bf Daley and Servi \cite{daley1992exploiting}:} Continues the track of the QIE, using taboo probabilities.
		
		{\bf Heyde \cite{heyde1992some}:} Quasi-likelihood estimation methods for stationary processes and queueing examples. 
		
		{\bf Jain \cite{jain1992relative}:}  Derives the relative efficiency of a parameter for the M/G/1 queueing system based on reduced and full likelihood functions. In addition, Monte Carlo simulations were carried out to study the finite sample properties for estimating the parameters of an M/G/1 queueing system. 
		
		{\bf Kumar \cite{kumar1992average}:} Studies the bias in the means of average idle time and average queue length estimates, over the interval [0, t], in a transient M/M/1 queue.

		{\bf Singpurwalla \cite{singpurwalla1992discussion}:} A discussion paper about \cite{thiruvaiyaru1992empirical}. The same issue for QUESTA also has a rejoinder for the discussion.
		
		{\bf Thiruvaiyaru and Basawa \cite{thiruvaiyaru1992empirical}:} Discusses empirical Bayes estimation for variations of M/M/1 queues and Jackson networks.\\

		\medskip
		{\bf  {\color{mycolor2}1993}}
		
		\medskip
		
		{\bf Daley and Servi \cite{daley1993two}:}  
		Discuss a fairly general Markov chain setting for
		describing a stochastic process at intermediate time points $r$ in $r \in (0, n)$ conditional
		on certain known behaviour of the process both on the interval and at the endpoints
		$0$ and $n$.

		{\bf Glynn, Melamed and Whitt \cite{glynn1993estimating}:} Constructs confidence intervals for estimators and perform statistical tests by establishing a joint central limit theorem for customer and time averages by applying a martingale central limit theorem in a Markov framework. \\
		
		\medskip
		{\bf  {\color{mycolor2}1994}}
		
		\medskip

		{\bf Armero and Bayarri \cite{armero1994prior}:} Presents a  Bayesian approach to  predict several quantities in an M/M/1 queue in equilibrium.
		
		{\bf Armero and Bayarri, M.J. \cite{armero1994bayesian}:} Bayesian ``prediction'' in M/M/1 queues is considered. The meaning is Bayesian inference for steady state quantities such as the distribution of queue lengths.
		
		{\bf Armero \cite{armero1994bayesianInf}:} Another Bayesian inference paper.
		
		\textbf{Chandra and Lee \cite{chandrs1994transactional}:} Presents  Bayesian methods for inferring customer behavior from transactional data in telecommunications systems. 
		
		{\bf Chen, Walrand and Messerschmitt \cite{chen1994parameter}:} Perhaps the first "probing" paper. Deals with arrivals in a deterministic service time queue and estimates the Poisson arrival rates based on probe delays.

		{\bf Jang and Liu \cite{jang1994waiting}:}  Presents a new queueing formula applicable in manufacturing which uses variables  easier to estimate than the variance such as the number of machine idle periods.
		
		{\bf Jones and Larson \cite{jones1994efficient}:} Develops an efficient algorithm for event probabilities of order statistics and uses it for the queue inference engine (\cite{larson1990queue}).
		
		{\bf Pitts \cite{pitts1994nonparametric}:} Analysis of non-parametric estimation of the GI/G/1 queue input distributions based on observation of the waiting time.\\
		
	
		\medskip
		{\bf  {\color{mycolor2}1995}}
		
		\medskip
		
		\textbf{Duffield, Lewis, O'Connell,  Russell,  and Toomey \cite{duffield1995entropy}:}
		Estimates directly the thermodynamic
		entropy of the data-stream at an input-port. From this,  the quality-of-service parameters can be calculated rapidly.
		
		{\bf Jain \cite{jain1995estimating}:} Change point detection in an M/M/1 queue.
		
		{\bf Muthu and Sampathkumar \cite{muthu1995estimation}:} The maximum likelihood estimates of the parameters involved in a finite capacity priority queueing model are obtained. The precision of the maximum likelihood estimates is studied using likelihood theory for Markov processes.
		
		{\bf Masuda \cite{masuda1995exploiting}:} Provides sufficient conditions under which the intuition (based on partial observations) can be justified, and investigates related properties of queueing systems.\\
		

		\medskip
		{\bf  {\color{mycolor2}1996}}
		
		\medskip
		
		{\bf Basawa,  Bhat, and Lund \cite{basawa1996maximum}:} MLE for GI/G/1 based on waiting time data.
		
		{\bf Dimitrijevic \cite{dimitrijevic1996inferring}:} Considers the problem of inferring the queue length of an M/G/1 queue using transactional data of a busy period.
		
		{\bf Manjunath and Molle \cite{manjunath1996passive}:} Introduces a new off-line estimation algorithm for the waiting times of departing customers in an M/G/1 queue with FCFS service by decoupling the arrival time constraints from the customer departure times.
		
		{\bf Massey, Parker, and  Whitt \cite{massey1996estimating}:
		} Estimate the parameters of a nonhomogeneous Poisson process with linear rate over a finite interval, based on the number of counts in measurement subintervals.

		{\bf Sohn \cite{sohn1996empirical}:} Simple M/M/1 Bayesian parameter estimation. 
		
		{\bf Sohn \cite{sohn1996influence}:} Bayesian estimation of M/M/1 using several competing methods. \\
		
		\medskip
		{\bf  {\color{mycolor2}1997}}
		
		\medskip

		{\bf Armero and Bayarri \cite{armero1997bayesian}:} Bayesian inference of M/M/$\infty$.
		
		{\bf Basawa,  Lund, and Bhat \cite{basawa1997estimating}:} Extends \cite{basawa1996maximum} using estimating functions.
		
		{\bf Bhat, Miller and Rao \cite{bhat1997statistical}:}  A survey paper, a decade after the previous Survey by Bhat and Rao, \cite{bhat1987statistical}.
		
		{\bf Daley and Servi \cite{daley1997estimating}:} Computes the distributions and moments of waiting times of customers within a busy period in an FCFS queuing system with a Poisson arrival process by exploiting an embedded Markov chain.

		{\bf Glynn and Torres \cite{glynn1997parametric}:} Deals with estimation of the tail properties of the workload process in both the M/M/1 queue and queues with more complex arrivals such as MMPP.

		{\bf Ho and Cassandras \cite{ho1997perturbation}:} A survey on perturbation theory.

		{\bf Pickands and Stine \cite{pickands1997estimation}:} Discrete time 
		M/G/$\infty$ queue.
		
		{\bf Toyoizumi \cite{toyoizumi1997sengupta}:} Waiting time inference in G/G/1 queues in a non-parametric manner using ``Sengupta's invariant relationship''.\\

		\medskip
		{\bf  {\color{mycolor2}1998}}
		
		\medskip
		{\bf Daley and Servi \cite{daley1998moment}:} Computes the distributions and moments of waiting times of customers within a busy period in an FCFS queuing system with a Poisson arrival process by exploiting an embedded Markov chain.
		
		{\bf Ganesh, Green, O'Connell and Pitts \cite{ganesh1998bayesian}:} Appears like a ``visionary'' paper on the use of non-parametric Bayesian methods in network management.
		
		{\bf Insua, Wiper and Ruggeri \cite{insua1998bayesian}:} Bayesian inference for M/G/1 queues with either Erlang or hyper-exponential service distributions.
		
		{\bf Mandelbaum and Zeltyn \cite{mandelbaum1998estimating}:} Queuing inference estimation in networks.

		{\bf Rodrigues and Leite \cite{rodrigues1998note}:} A Bayesian inference about the traffic intensity  in an M/M/1 queue, without worrying about nuisance parameters. 
		
		{\bf Sharma and Mazumdar \cite{sharma1998estimating}:}  Proposes several schemes that the call acceptance controller, at the entering node of an ATM network, can use to estimate the traffic of the users on the various routes in the network by sending a probing stream.
		
		{\bf Wiper \cite{wiper1998bayesian}:} Perhaps complements \cite{insua1998bayesian} with analysis of G/M/c queues with the G being Erlang or hyper-exponential renewal processes. \\
		
		\medskip
		{\bf  {\color{mycolor2}1999}}
		
		\medskip
		
		{\bf Acharya \cite{acharya1999normal}:} Analyses the rate of convergence of the distribution of MLEs in GI/G/1 queues with assumptions on the distributions as being from exponential families.

		{\bf Conti \cite{conti1999large}:} Baysian inference for a Geo/G/1 Discrete time queue.
		
		{\bf Bingham and Pitts \cite{bingham1999non}:} Non-parameteric estimation in 
		M/G/$\infty$ queues.
		
		{\bf Bingham and Pitts \cite{bingham1999nonparametric}:}  Estimates the arrival rate of an M/G/1 queue given observations of the busy and idle periods of this queue.

		{\bf Jones \cite{jones1999inferring}:}  Analyses queues in the presence of balking, using only the service start and stop data utilized in Larson's Queue Inference Engine.

		{\bf Rodrigo and Vazquez \cite{rodrigo1999large}:} Analyses a general G/G/1 retrial queueing systems from a statistical viewpoint.
		
		{\bf Sharma \cite{sharma1999estimatingB}:} Using the measurement tools available on the Internet, suggests and compares different
		estimators to estimate aggregate traffic intensities at
		various nodes in the network. \\

		\medskip
		{\bf  {\color{mycolor2}2000}}
		
		\medskip
		
		{\bf Armero and Conesa \cite{armero2000prediction}:} Statistical analysis of bulk arrival queues from a Bayesian point of view.
		
		{\bf Duffield \cite{duffield2000large}:} 
		Analyses the impact of measurement error within the framework of
		Large Deviation theory.
		
		{\bf Glynn and Zeevi \cite{glynn28estimating}:} Estimates tail probabilities in queues.

		{\bf Jain \cite{jain2000autoregressive}:} Sequential analysis.

		{\bf Jain and Rao \cite{jain2000statistical}:} Investigates the problems of statistical inference for the GI/G/1 queueing system. 
		
		{\bf Zheng and  Seila \cite{zheng2000some}:} 
		Construct  estimators for the  limiting expected number of customers in the queue (and several other performance measures)  with better sampling properties in comparison to the existing estimators.\\
		
		\medskip
		{\bf  {\color{mycolor2}2001}}
		
		\medskip
		
		{\bf Alouf, Nain and Towsley \cite{alouf2001inferring}:} Probing estimation for M/M/1/K queues using moment-based estimators based on a variety of computable performance measures. 
		
		{\bf Huang and Brill \cite{huang2001estimation}:} Deriving the minimum variance unbiased estimator (MVUE) and the maximum likelihood estimator (MLE) of the stationary probability function of the number of customers in a collection of independent M/G/c/c subsystems. 
		
		{\bf Jang, Suh and Liu \cite{jang2001new}:} Presents a new GI/G/2 queueing formula which uses a slightly different set of data easier to obtain than the variance of inter-arrival time. 
		
		{\bf Paschalidis and Vassilaras \cite{paschalidis2001estimation}:} Buffer overflow probabilities in queues with correlated arrival and service processes using large deviations. \\

		
		\medskip
		{\bf  {\color{mycolor2}2002}}
		
		\medskip

		{\bf Conti \cite{conti2002nonparametric}:} Non-parametric statistical analysis of a discrete-time queueing system is considered and estimation of performance measures of the system is studied. 
		
		{\bf Conti and De Giovanni \cite{conti2002queueing}:}  Considers performance evaluation of a discrete-time GI/G/1  queueing model with a focus on the equilibrium distribution of the waiting time. 
		
		{\bf Sohn \cite{sohn2002robust}:} Even though the title has ``Robust''', this paper appears to be a standard M/M/1 Baysian inference paper using the input data.
		
		{\bf Zhang, Xia, Squillante and Mills \cite{zhang2002workload}:} A  general approach to infer the per-class service times at different servers in multi-class queueing models.\\
		
		\medskip
		{\bf  {\color{mycolor2}2003}}
		
		\medskip
		
		{\bf Cao, Andersson, Nyberg, and Kihl \cite{cao2003web}:}
		Web server performance is modeled via an M/G/1/K*PS queue for which the authors also carry out maximum likelihood estimation of the parameters.
		
		
		{\bf Pichitlamken, Deslauriers, L'Ecuyer, and Avramidis \cite{pichitlamken2003modelling}:} This is a simulation modelling paper where the authors also carry out parameter estimation for the simulation model data via a real data set. 
		\\
		
		\medskip
		{\bf  {\color{mycolor2}2004}}
		
		\medskip
		
		{\bf Aus{\'\i}n, Wiper and Lillo \cite{ausin2004bayesian}:} Bayesian inference of M/G/1 using phase type representations of the G.
		
		{\bf Conti \cite{conti2004bootstrap}:} A Bayesian non-parametric approach to the analysis of discrete-time queueing models. 
		
		{\bf Fearnhead \cite{fearnhead2004filtering}:} Using forward-backward algorithm to do inference for  M/G/1 and Er/G/1 queues.

		{\bf Hall and Park \cite{hall2004nonparametric}:} An M/G/$\infty$ non-parametric paper.
		
		{\bf Wang,  Chen, and  Ke \cite{estimates2004kuo}:} Maximum likelihood estimates and confidence intervals of an M/M/R/N queue with balking and heterogeneous servers.\\

		\medskip
		{\bf  {\color{mycolor2}2005}}
		
		\medskip
		
		{\bf Bladt and S{\o}rensen \cite{bladt2005statistical}:}  Likelihood inference for discretely observed Markov jump processes with finite state space. 
		
		{\bf Brown, Gans, Mandelbaum, Sakov, Shen, Zeltyn and Zhao \cite{brown2005statistical}:} Major paper dealing with telephone call centre data analysis.
		
		{\bf Hei, Bensaou and Tsang \cite{hei2005light}:} Probing focusing on the inter-departure SCV of the probing stream in tandem finite buffer queues.

		{\bf Mandjes and van de Meent \cite{mandjes2005inferring}:}
		Propose an approach to accurately determine the burstiness of a network link on small time-scales (for instance 10 ms), by sampling the buffer occupancy (for instance) every second.
		
		{\bf Prieger \cite{prieger2005estimation}:} Shows that the MLE based on the complete inter-arrival and service times (IST) dominates the MLE based on the number of units in service (NIS), in terms of ease of implementation, bias, and variance.
		
		{\bf Ross and Shanthikumar \cite{ross2005estimating}:} Estimating effective capacity in Erlang loss systems
		under competition.
		
		{\bf Neuts \cite{neuts2005reflections}:} Reflections on statistical methods for complex stochastic systems.\\
		
		\medskip
		{\bf  {\color{mycolor2}2006}}
		
		\medskip
		
		{\bf Castellanos, Morales, Mayoral, Fried and Armero  \cite{castellanos2006bayesian}:} Develops a Bayesian analysis of queueing systems in applications of the machine interference problem, like job-shop type systems, telecommunication traffic, semiconductor manufacturing or transport.

		{\bf Chick \cite{chick2006subjective}:} A survey chapter on subjective probability and the Bayesian approach, specifically in Monte-Carlo simulation, yet gives some insight into queueing inference.

		{\bf Chu and Ke \cite{chu2006confidence}:} Construction of confidence intervals of mean response time for an M/G/1 FCFS queueing system.

		{\bf Doucet, Montesano Jasra \cite{doucet2006optimal}:}  Presents a  trans-dimensional Sequential Monte Carlo method for online Bayesian inference in partially observed point processes.

		{\bf Hansen and Pitts \cite{hansen2006nonparametric}:} Non-parametric estimation of the service time distribution and the traffic intensity in M/G/1 queues based on observations of the workload.
		
		{\bf Hei, Bensaou and Tsang \cite{hei2006model}:}  Similar to \cite{hei2005light} but here the focus is on inter-departure SCV of the two consecutive probing packets.

		{\bf Ke and Chu \cite{ke2006nonparametric}:}  Proposes a consistent and asymptotically normal estimator of intensity for a queueing system with distribution-free inter-arrival and service times. 
		
		{\bf Kim, Shiravi, and Min \cite{kim2006congestion}:} This paper deals with network congestion analysis via approximating processes and uses real network trace data for parameter fitting of self-similar network traffic using the index of dispersion for counts and coefficient of determination.
		
		{\bf Liu, Heo, Sha and Zhu \cite{liu2006adaptive}:} Proposes a queueing-model-based adaptive control approach for controlling the performance of computing systems.
		
		{\bf Liu, Wynter, Xia and Zhang \cite{liu2006parameter}:} presents an approach for solving the problem of calibration of model parameters in the queueing network framework using inference techniques.
		
		{\bf Rodrigo \cite{rodrigo2006estimators}:} Analyse the M/G/1 retrial queue from a statistical viewpoint.

		{\bf Wang, Ke, Wang and Ho \cite{wang2006maximum}:}
		Studies MLE and confidence intervals of an M/M/R queue with heterogeneous servers under steady-state conditions.\\

		\medskip
		{\bf  {\color{mycolor2}2007}}
		
		\medskip
		
		{\bf Aus{\'\i}n, Lillo and Wiper \cite{ausin2007bayesian}:} Considers the problem of designing a GI/M/c queueing system.

		{\bf Chu and Ke  \cite{chu2007interval}:} Proposes a consistent and asymptotically normal  estimator of the mean response time for a G/M/1 queueing system, which is based on the empirical Laplace function. 
		
		{\bf Chu and Ke \cite{chu2007mean}:} Estimation and confidence
		interval of mean response time for a G/G/1 queueing system using data-based recursion relation and bootstrap methods.

		{\bf Morales, Castellanos, Mayoral, Fried and Armero \cite{morales2007bayesian}:}  Exploits Bayesian criteria for designing an M/M/c//r queueing system with spares.

		{\bf Park \cite{park2007choice}:}  The use of auxiliary functions in non-parametric inference for the M/G/$\infty$ queueing model is considered. 
		
		{\bf Ross, Taimre and Pollett \cite{ross2007estimation}:} Estimation of rates in M/M/c queues using observations at discrete queues and MLE estimates of an approximate Orenstein Ullenbeck (OU) process.\\
		

		\medskip
		{\bf  {\color{mycolor2}2008}}
		
		\medskip
		
		{\bf Aus{\'\i}n,  Wiper and Lillo \cite{ausin2008bayesian}:} Bayesian inference for the transient behaviour and duration of a busy period in a single server queueing
		system with general, unknown distributions for the inter-arrival and service times is investigated. 
		
		{\bf Basawa, Bhat and Zhou \cite{basawa2008parameter}:} Parameter estimation based on the differences of two positive exponential family random variables is studied.
		
		{\bf Casale, Cremonesi and Turrin \cite{casale2008robust}:} Proposed service
		time estimation techniques based on robust and constrained optimization.
		
		\textbf{Casale, Zhang and  Smirni \cite{casale2008interarrival}:}
		Several contributions to the Markovian traffic analysis. 
		
		{\bf Choudhury and Borthakur \cite{choudhury2008bayesian}:} Bayesian-based techniques for analysis of the M/M/1 queueing model.
		
		{\bf Dey \cite{dey2008note}:} Bayes' estimators of the traffic intensity and various queue characteristics in an M/M/1 queue under quadratic error loss function  have been derived. 
		
		{\bf Ke, Ko and Sheu \cite{ke2008estimation}:} Proposes an estimator for the expected busy period  of a controllable
		M/G/1 queueing system in which the server applies a bicriterion policy during his idle
		period.

		{\bf Ke, Ko and Chiou \cite{ke2008analysis}:}
		Presents a sensitivity investigation of the expected busy period
		for a controllable M/G/1 queueing system by means of a factorial design statistical
		analysis. 
		
		{\bf Kim and Park \cite{kim2008new}:} Introduces methods of queue
		inference which can find the internal behaviours of queueing
		systems with only external observations, arrival and departure time. 
		
		{\bf Sutton and Jordan \cite{sutton2008probabilistic}:} Analysing queueing networks from the probabilistic modelling perspective, applying inference methods from graphical models that afford significantly more modelling flexibility.
		
		{\bf Ramirez, Lillo and Wiper \cite{ramirez2008bayesian}:}
		Considers a mixture of two-parameter Pareto distributions as a
		model for heavy-tailed data and use a Bayesian approach based on the birth-death
		Markov chain Monte Carlo algorithm to fit this model.\\
		
		\medskip
		
		{\bf  {\color{mycolor2}2009}}
		\medskip

		{\bf Baccelli, Kauffmann and Veitch \cite{baccelli2009inverse}:} Points out the importance of inverse problems in queueing theory, which aim
		to deduce unknown parameters of the system based on partially observed trajectories.

		{\bf Baccelli, Kauffmann and Veitch \cite{baccelli2009towards}:} Evaluates the algorithm proposed in \cite{baccelli2009inverse}.
		
		\textbf{Comert and Cetin \cite{comert2009queue}} Presents a real-time estimation of queue lengths
		from the location information of probe vehicles in a queue at an isolated and under-saturated intersection.
		
		{\bf Chu and Ke \cite{chu2009analysis}:} Constructs confidence intervals of intensity for a queueing system, which are based on four different bootstrap methods.

		{\bf Duffy and Meyn \cite{duffy2009estimating}:}   Conjectures  and presents support for this: a consistent sequence of non-parametric estimators can be constructed that satisfies a large deviation principle. 
		
		{\bf Gorst-Rasmussen Hansen \cite{gorst2009asymptotic}:} Proposes a framework based on empirical process techniques
		for inference about waiting time and patience distributions in multi-server queues with abandonment.
		
		{\bf Heckmuller and Wolfinger \cite{heckmuller2009reconstructing}:} Proposes methods to estimate the
		parameters of arrival processes to G/D/1 queueing systems only based on observed departures from the system.

		{\bf Ibrahim and Whitt \cite{ibrahim2009real}:}
		Studies the performance of alternative real-time
		delay estimators based on recent customer delay experiences.

		{\bf Kiessler and Lund \cite{kiessler2009technical}:}  A note that considers traffic intensity estimation in the classical M/G/1 queue.

		{\bf Ke and Chu \cite{ke2009comparison}:} Proposes a consistent and asymptotically normal
		estimator of intensity for a queueing system with distribution-free inter-arrival
		and service times.

		{\bf Kraft, Pacheco-Sanchez, Casale and Dawson \cite{kraft2009estimating}:} Proposes a linear regression method and a maximum likelihood
		technique for estimating the service demands of requests based on the measurement of their response times instead of their CPU utilization. 
		
		\textbf{Liu,  Wu, Ma, and Hu \cite{liu2009real}:}
		Presents an approach to estimate time-dependent queue length even when the signal links are congested.
		
		{\bf Mandjes and {\.Z}uraniewski \cite{mandjes2009queueing}:} Develops  queueing-based procedures to  (statistically) detect overload in communication networks, in a setting in which each connection consumes roughly the same amount of bandwidth.
		
		{\bf Mandjes and van De Meent \cite{mandjes2009resource}:} The focus is on dimensioning as the approach for delivering
		performance requirements of the network. 
		
		{\bf Nam, Kim and Sung \cite{nam2009estimation}:}
		Estimates the available bandwidth for an M/G/1 queueing system.

		{\bf Novak and Watson \cite{novak2009determining}:} Presents a technique to estimate the arrival rate from delay measurements, acquired using single-packet probing. \\

		\medskip
		{\bf  {\color{mycolor2}2010}}
		
		\medskip

		{\bf Chen and Zhou \cite{chen2010simulation}:}
		Proposes a non-linear quantile regression
		model for the relationship between stationary cycle time quantiles and corresponding throughput rates of a manufacturing system.
		
		{\bf Duffy and Meyn \cite{duffy2010most}:}  Deals with large deviations showing that in broad generality, that estimates of the steady-state mean position of a reflected random walk have a high likelihood of over-estimation.

		{\bf Frey and Kaplan \cite{frey2010queue}:} Introduces an algorithm for queue inference problems involving periodic reporting data.
		
		{\bf Pin, Veitch, and Kauffmann \cite{pin2010statistical}:} This paper incorporates queueing theory results in the application of network tomography where statistical estimation of delays in a multicast tree using an EM algorithm is developed.
		
		\textbf{Gans, Liu, Mandelbaum, Shen, and Ye \cite{gans2010service}:} Studies operational heterogeneity of call center agents where the
		proxy for heterogeneity is agents' service times (call durations).
		
		{\bf Heckmueller and Wolfinger \cite{heckmueller2010reconstructing}:} Proposes methods to estimate the parameters of arrival processes to G/D/1 queueing systems only based
		on observed departures from the system.
		
		{\bf Pin, Veitch and Kauffmann \cite{pin2010statistical}:} Focuses on a specific delay tomographic problem on a multicast diffusion
		tree, where end-to-end delays are observed at every leaf of the tree, and mean
		sojourn times are estimated for every node in the tree.
		
		{\bf Ramirez-Cobo, Lillo, Wilson and Wiper \cite{ramirez2010bayesian}:} Presents a method for carrying out Bayesian estimation
		for the double Pareto lognormal  distribution.

		{\bf Sutton and Jordan \cite{sutton2010learning}:}
		Presents a viewpoint that combines queueing networks and graphical models, allowing Markov chain Monte Carlo to be applied.

		{\bf Xu, Zhang and Ding \cite{xu2010hypothesis}:}
		Discusses testing hypotheses and confidence regions with correct levels
		for the mean sojourn time of an M/M/1 queueing system.
		
		{\bf Zhang and Xu \cite{zhang2010confidence}:} Discuss constructing confidence intervals  of performance measures for an M/G/1 queueing system.

		{\bf Zuraniewski, Mandjes and Mellia \cite{zuraniewski2010empirical}:} Explores techniques for detecting unanticipated load changes with a focus on large-deviations based techniques developed earlier in \cite{mandjes2009queueing}.\\

		\newpage
		\medskip
		{\bf  {\color{mycolor2}2011}}
		
		\medskip
		
		{\bf Abramov \cite{abramovstatistical}:} Statistical bounds for certain output characteristics of the M/GI/1/n and GI/M/1/n loss queueing systems are derived on the basis of large samples of an input characteristic of these systems.
		
		{\bf Amani, Kihl and Robertsson \cite{amani2011multi}:} An applications paper to computer systems.

		\textbf{Ban, Hao, and Sun \cite{ban2011real}:}
		Studies how to estimate real time queue lengths at signalized intersections using the intersection
		travel times collected from mobile traffic sensors.

		{\bf Chen, Nan, Zhou \cite{chen2011performance}:}
		Investigates the statistical process control application for monitoring queue length data in M/G/1 systems.

		{\bf Feng, Dube, and Zhang \cite{feng2011estimation}:} 
		Considers estimation problems in G/G/$\infty$ queue under
		incomplete information. Specifically, where it is infeasible to track each individual job
		in the system and only aggregate statistics are known or observable.

		{\bf Veeger,  Etman,  Lefeber, Adan,   Van Herk,  and Rooda \cite{veeger2011predicting}:}
		Predicting cycle time distributions for integrated processing workstations: an aggregate modelling approach.

		{\bf Gr{\"u}bel and Wegener \cite{grubel2011matchmaking}:}  In M/G/$\infty$ systems, considers the matching and exponentiality problems where the only observations are the order statistics associated with the departure
		times  and the order in which the customers arrive and depart, respectively.

		{\bf Ibrahim and Whitt \cite{ibrahim2011wait}:}
		Develops real-time delay predictors for many-server service systems with a time-varying arrival rate, a time-varying number of servers, and customer abandonment.

		{\bf Mandjes and Zuraniewski \cite{mandjes2011m}:} M/G/$\infty$ change point detection using large deviations. 
		
		{\bf Manoharan and Jose \cite{manoharan2011markovian}:}  Considers an M/M/1 queueing system with customer impatience in the form of random balking.
		
		{\bf McCabe, Martin and Harris \cite{mccabe2011efficient}:}  Presents an efficient probabilistic forecast of integer-valued random variables that can be interpreted
		as a queue, stock, birth and death process or branching process.
		
		{\bf Park, Kim and Willemain \cite{park2011analysis}:} Proposes  new approaches that can analyse the unobservable
		queues using external observations.
		
		{\bf Sousa-Vieira \cite{sousa2011suitability}:} Considers the suitability of the M/G/$\infty$ process for modelling the spatial and quality scalability extensions of the H.264 standard in video traffic modelling.
		
		{\bf Srinivas, Rao and Kale  \cite{srinivas2011estimation}:} Maximum likelihood and uniform minimum variance unbiased estimators of measures in the M/M/1 queue are obtained and compared.

		{\bf Sutton and Jordan \cite{sutton2011bayesian}:} A Bayesian inference paper by computer systems researchers. \\

		\medskip
		{\bf  {\color{mycolor2}2012}}

		\medskip
		
		{\bf Duffy and Meyn \cite{duffy2012large}:} Large deviation asymptotics for busy periods for a queue.
		
		{\bf  Fabris-Rotelli, Kraamwinkel,  and et al \cite{fabris2012investigation}:} A historical and theoretical overview of  G/M and M/G
		queueing processes.
		
		{\bf Hu and Lee \cite{hu2012parameter}:} Consider a parameter estimation problem when the state process is a reflected fractional
		Brownian motion (RFBM) with a non-zero drift parameter and the observation is the associated
		local time process. 
		
		{\bf Jones \cite{jones2012remarks}:} Remarks on queue inference from departure data alone and the importance of the queue inference engine.
		
		{\bf Kauffmann \cite{kauffmann2012inverse}:}
		Proposes a new approach, on the basis of existing TCP connections and reaching therefore a zero probing overhead based on the theory of inverse
		problems in bandwidth sharing networks.

		{\bf Kim and Whitt \cite{kim2012statistical}:} 
		Statistical analysis with Little's law.
		
		{\bf Kim and Whitt \cite{kim2012statisticalSupplementary}:} Contains supplementary martial to  \cite{kim2012statistical}.
		
		{\bf Kim and Whitt \cite{kim2012estimating}:}
		Estimating waiting times with the time-varying Little's law.
		
		{\bf Kim and Whitt \cite{kim2012estimatingAppendix}:} Appendix to  \cite{kim2012estimating}.

		{\bf McVinish and Pollett  \cite{mcvinishconstructing}:} Uses estimating equations to get estimators for M/M/c queues and related models. Performance is compared to  \cite{ross2007estimation}.
		
		{\bf Mohammadi and Salehi-Rad \cite{mohammadi2012bayesian}:} Exploits the Bayesian inference and prediction for an M/G/1 queuing model with optional second re-service.

		{\bf Nelgabats, Nov and Weiss \cite{nelgabats2012}:} M/G/$\infty$ estimation.
		
		{\bf Ren and Li \cite{ren2012bayesB}:} Bayesian estimator of the traffic intensity in an M/M/1 queue is derived under a new weighted square error loss function. 
		
		{\bf Ren and Wang \cite{ren2012bayes}:}  Similer to \cite{ren2012bayesB}. Bayesian estimators of the traffic intensity in an M/M/1 queue are derived under a precautionary loss function. 
		
		{\bf Whitt \cite{whitt2012fitting}:} Fitting birth-and-death queueing model to data.
		\\
		
		\medskip
		{\bf  {\color{mycolor2}2013 }}
		\medskip
		
		{\bf Larson \cite{larson2013queue}:} A brief review on QIE.

		{\bf Acharya,  Rodríguez-Sánchez and  Villarreal-Rodríguez \cite{acharya2013maximum}:}
		Presents the derivation of maximum likelihood estimates for the arrival rate and service rates in a stationary M/M/c queue with heterogeneous servers.

		\textbf{Chow \cite{chow2013observable}:}
		Analysis of queueing model based on chaotic mapping.

		\textbf{Li, Chen, Li and Zhang \cite{li2013freeway}:} Proposes a
		new algorithm based on the temporal--spatial queueing model to describe
		the fast travel-time variations using only the speed and headway time
		series that is measured at upstream and downstream detectors.
		
		{ \bf Wiler, Bolandifar, Griffey, Poirier, and Olsen  \cite{wiler2013emergency}:} This paper applies queueing theory and inference of model parameters to derive and validate a novel queuing theory-based model that predicts
		the effect of various patient crowding scenarios on patients.
		
		\textbf{Weerasinghe and Mandelbaum \cite{weerasinghe2013abandonment}:} studies the tradeoff between blocking and abandonment, with
		cost accumulated over a random, finite time horizon of a controlled queueing system of the G/M/n/B+M type with many
		servers and impatient customers. 
		\\
		
		\medskip
		{\bf  {\color{mycolor2}2014 }}
		
		\medskip
		
		{\bf Antunes, Jacinto and Pacheco \cite{antunes2014probing}:}
		estimation of the arrival rate and the service
		time moments of an M/G/1 queue with probing, i.e.,
		with special customers (probes) entering the system. The
		probe inter-arrival times are i.i.d. and probe service times
		follow a general positive distribution. The only observations
		used are the arrival times, service times and departure
		times of probes. We derive the main equations from which
		the quantities of interest can be estimated. Two particular
		probe arrivals, deterministic and Poisson, are investigated.

		\textbf{Azriel, Feigin and Mandelbaum \cite{azriel2014erlang}:}
		Proposes a new model called Erlang-S, where ``S'' stands for Servers where there is a pool of present servers, some
		of whom are available to serve customers from the queue while others are not, and the process
		of becoming available or unavailable is modelled explicitly.

		{\textbf{Dinha, Andrewa and Nazarathy \cite{dinh2014architecture}:}
			A conceptual and numerical contribution to design
			and control of speed-scaled systems in view of parameter uncertainty.
			
			{\bf Edelmann and Wichelhaus \cite{edelmann2014nonparametric}:}
			Two different nonparametric estimation approach  for discrete-time stochastic networks of
			Geom$^X$/G/$\infty$ queues where the observation consists of the external arrival and external departure processes at the nodes over some time.

			{\bf He, Li, Huang and Lei \cite{he2014maximum}:}  Considers
			the queuing system as a black box and derive a performance index for the queuing system by the principle of maximum entropy only
			on the assumption that the queue is stable.
			
			{\bf Kannan and Jabarali \cite{kannan2014parameter}:}
			This paper deals with maximum likelihood estimation parameters for a variant of an M/M/1 queue with vacations. 
			
			{\bf Kim and Whitt \cite{kim2014choosing}:}
			Considers different issues in testing the suitability of the nonhomogeneous  Poisson process as an arrival process with service system data. 
			
			{\bf Kim and Whitt \cite{kim2014call}:} (in the following of \cite{kim2014choosing})
			shows that call center and hospital arrivals are well modelled by nonhomogeneous Poisson processes.
			
			{\textbf{Senderovich,  Weidlich,  Gal, and  Mandelbaum \cite{senderovich2014queue}:}
				Establish a queueing perspective in operational process mining and demonstrate the value of queue mining using the specific operational problem of online delay prediction.

				\textbf{Yom-Tov and Mandelbaum \cite{yom2014erlang}:} Analyses a queueing model, where customers can return to service several times during their sojourn within the system.
				\\
				
				\medskip
				{\bf  {\color{mycolor2}2015} }

				\medskip
				
				\textbf{Bakholdina1 and Gortsev \cite{bakholdina2015optimal}:} Focused on the problem of optimal estimation of the
				states of the modulated
				semi-synchronous integrated flow of events.
				
				\textbf{Burkatovskaya, Kabanova,  and Vorobeychikov  \cite{vorobeychikov2015cusum}:}
				CUSUM algorithms for parameter estimation
				in queueing systems where the arrival process is a Markov-modulated Poisson process.
				
				\textbf{Cahoy, Polito, and Phoha \cite{cahoy2015transient}:} Statistical analysis of fractional M/M/1 queue and fractional birth-death processes; the point processes governed by difference differential
				equations containing fractional derivative operators. 
				
				{\bf Chen and Zhou \cite{chen2015cusum}:} Propose the
				cumulative sum (CUSUM) schemes to efficiently monitor the performance of typical queueing systems based on different sampling schemes.
				
				\textbf{Dong and Whitt \cite{dong2015stochastic}:}
				Explores a stochastic grey-box modelling of queueing systems by fitting birth-and-death processes to data.
				
				\textbf{Dong and Whitt \cite{dong2015using}:}
				Using a birth-and-death process to estimate the steady-State distribution
				of a periodic queue.
				
				\textbf{Efrosinin, Winkler, and Martin \cite{efrosinin2015confidence}:} Considers the problem of estimation and confidence interval construction of a Markovian
				controllable queueing system with unreliable server and constant retrial policy.

				\textbf{Goldenshluger \cite{goldenshluger2015nonparametric}:} Non-parametric estimation of service time distribution of the M/G/1 queue from incomplete data on the queue. 
				
				\textbf{Gurvich, Huang and Mandelbaum \cite{gurvich2013excursion}:} Proposes a diffusion approximation for a  many-server Erlang-A queue.
				
				{\bf Kim and Whitt \cite{kim2015power}:} Similar to \cite{kim2014choosing}   
				
				\textbf{Liu, Wu, and Michalopoulos  \cite{liu2015improving}:} Improves queue size estimation by  proposing different ramp queue estimation algorithms.
				
				\textbf{Mohajerzadeh, Yaghmaee, and Zahmatkesh \cite{mohajerzadeh2015efficient}:} Proposed a method to  prolong the network lifetime and to estimate the target parameter efficiently in wireless sensor networks.
				
				\textbf{Schweer and Wichelhaus \cite{schweer2015nonparametric}} Estimation of the service time distribution in the discrete-time GI/G/$\infty$-queue based solely on information on the arrival and departure processes is considered. The focus is put on the estimation approach via the so called``sequence of difference'' and proving a functional central limit theorem for the resultant estimator.

				\textbf{Senderovich, Leemans, Harel,  Gal, Mandelbaum, and van der Aalst \cite{senderovich2015discovering}:}  Explores the influence of available information in the log 
				on the accuracy of the queue mining techniques.

				\textbf{Senderovich, Weidlich, Gal, Mandelbaum \cite{senderovich2015queue}:} Queue mining for delay prediction in multi-class
				service processes.
				
				\textbf{Srinivas and Kale \cite{srinivas2015ml}:} 
				Compares the Maximum Likelihood (ML) and Uniformly Minimum Variance Unbiased (UMVU) estimation for the  M/D/1 queueing systems.

				\textbf{Sutartoa and Joelianto \cite{sutarto2015modeling}:} Presents an overview of urban traffic flow from the perspective of system theory and stochastic control.

				\textbf{Wang and Casale \cite{wang2015maximum}:} Proposes maximum likelihood (ML) estimators for service demands in closed queueing networks with load-independent and load-dependent stations. 
				
				\textbf{Wang, Pérez, and Casale \cite{wang2015filling}:} A software for parameter estimation.
				
				\textbf{Whitt \cite{whitt2015many}:}  Sequel to \cite{dong2015stochastic} and \cite{doucet2006optimal}. 
				Establishes many-server heavy-traffic fluid limits for the steady-state distribution
				and the fitted birth and death rates in periodic Mt/GI/$infty$ models.
				
				\textbf{Zhan, Li, and Ukkusuri \cite{zhan2015lane}:}
				In the context of transportation engineering, with the complete arrival and departure
				information, a car-following based simulation scheme is applied to estimate the
				real-time queue length for each lane.
				\\
				
				\medskip
				{\bf  {\color{mycolor2}2016} }
				\medskip
				
				\textbf{Amini,  Pedarsani,  Skabardonis, and Varaiya \cite{amini2016queue}:} Queue-length estimation using real-time traffic data.

				\textbf{Antunes, Jacinto, Pacheco, and Wichelhaus \cite{antunes2016estimation}:}
				Uses a probing strategy to estimate the time dependent
				traffic intensity in an M$_t$/G$_t$/1 queue, where the arrival rate and the general service-time distribution change from one time interval to another, and derive statistical properties of
				the proposed estimator. 
				
				\textbf{Anusha,  Sharma, Vanajakshi, Subramanian, and Rilett \cite{anusha2016model}:}
				Develops a model-based scheme to estimate the number of vehicles in queue and the total delay.
				
				\textbf{Comert \cite{comert2016queue}:} This paper, motivated by the field of traffic engineering, develops estimators for market penetration level and arrival rates based on  queue lengths from probe vehicles at isolated traffic intersections.
				
				\textbf{Cruz,  Quinino and Ho \cite{cruz2016bayesian}:}
				Uses a Bayesian technique, the sampling/importance resampling method  to estimate the parameters of  multi-server queueing systems in which inter-arrival and service times are exponentially distributed.
				
				\textbf{Ghorbani-Mandolakani and Salehi Rad \cite{ghorbani2016ml}:}
				Derives the ML and Bayes estimators of traffic intensity and asymptotic confidence intervals for mean system size of a two-phase tandem queueing
				model with a second optional service and random feedback and two heterogeneous servers.
				
				\textbf{Krishnasamy, Sen,   Johari,  and Shakkottai \cite{krishnasamy2016regret}:} Considers regret analysis of a server allocation problem where service rates of servers are unknown. This analysis is in the context of multi-armed bandits.
				
				\textbf{Morozov,  Nekrasova, Peshkova, and Rumyantsev \cite{morozov2016regeneration}:} Develops a novel approach to confidence estimation of the stationary measures in high performance multi-server queueing systems. 
				
				\textbf{Quinino and Cruz \cite{quinino2016bayesian}:} Describes a  Bayesian method for sample size determination for traffic intensity estimation.
				
				\textbf{Zammit, Fabri and Scerri \cite{zammit2016joint}:} A self-estimation algorithm is presented to
				jointly estimate the states and model
				parameters.
				
				\textbf{Zhang, Xu, and Mi \cite{zhang2016generalized}:} Considers the hypothesis tests of performance measures for an M/E$_k$/1 queueing system. \\
				
				\medskip
				{\bf  {\color{mycolor2}2017} }
				\medskip

				{\bf Cruz, Quinino,  and  Ho \cite{cruz2017bayesian}:}
				Bayesian estimation of traffic intensity based on queue length in a multi-server $M/M/s$ queue.

				{\bf den Boer and  Mandjes \cite{den2017convergence}:} Considers convergence rates of Laplace-transform based estimators.
				
				{\bf Gu, Qian,  and  Zhang \cite{gu2017traffic}:} A Bayesian probabilistic model along with an expectation–maximization extended Kalman filter (EM-EKF) algorithm is proposed for  traffic state estimation of urban road networks using a link queue model.
				
				{\bf Kim, Whitt, and Cha \cite{kim2017data}:}
				A data-driven model of an appointment-generated
				arrival process at an outpatient clinic.
				
				{\bf  Li, Tang, Yao, and Li \cite{doi:10.3141/2623-06}:}
				Proposes a cycle-by-cycle queue length estimation method using only probe data without the foregoing assumption for signalized intersections.

				{\bf  Quinino and Cruz \cite{quinino2017bayesian}:} A Bayesian method is described for sample size determination for traffic intensity estimation of an M/M/1 queue.

				{\bf Sutarto, Joelianto,  and Nugroho \cite{sutarto2017developing}:}
				Developing a stochastic model of queue length at a signalized
				intersection.
				
				{\bf Whitt and Zhang \cite{whitt2017data}:}
				Develops an aggregate stochastic model of an emergency department based on a careful study of data on individual patient arrival times
				and length of stay.
				
				{\bf Zammit, Fabri, and Scerri \cite{8317685}:} 
				Online state and multidimensional parameter estimation for a macroscopic model of a traffic junction based on the Expectation-Maximization algorithm and multidimensional Robbins-Monro stochastic approximation.
				\\
				
				\medskip
				{\bf  {\color{mycolor2}2018} }
				\medskip

				{\bf An, Wub, Xiaa, and Huanga \cite{an2018real}:} This paper from the field of traffic engineering focuses on real-time queue length estimation in the context of signalized intersections. 
				
				{\bf Almeida and Cruz \cite{almeida2018note}:} The Jeffreys prior is proposed to obtain the posterior and predictive distributions of some parameters of interest. Samples are obtained through simulation and some performance characteristics are analyzed.
				
				{\bf Cruz,  Almeida,  D’Angelo,  and van Woensel \cite{cruz2018traffic}:}
				Investigating the finite-sample behaviour of some well-known methods for the estimation of single-server finite Markovian queues or, in Kendall notation,  queues, namely, the maximum likelihood estimator, Bayesian methods, and bootstrap corrections.
				
				
				{\bf  Ozawa \cite{ozawa2019stability}:} Analysis of the stability condition of a two-dimensional QBD process and its
				application to evaluate the efficiency of  two-queue models.
				
				{\bf Polson and Sokolov \cite{7850972}:}
				Develops an efficient particle learning algorithm for real time online inference of states and parameters. This requires a two-step approach, first resampling the current particles with a mixture predictive distribution and second propagation of states using the conditional posterior distribution.
				
				{\bf Suyama,  Quinino,  and Cruz \cite{suyama2018simple}:} Estimators for the parameters of the Markovian multiserver queues are presented, from samples that are the number of clients in the system at arbitrary points and their sojourn times.\\
				
				
				
				\medskip
				{\bf  {\color{mycolor2}2019} }
				\medskip
				
				
				{\bf Bhat and Basawa\cite{bhat2019maximum}:} An overview of the literature on the use
				of the maximum likelihood method for estimating parameters in
				queueing models. 
				
				{\bf Emami,  Sarvi,  and Asadi Bagloee \cite{emami2019neural}:} A neural network algorithm for queue length estimation based on the concept of k-leader connected vehicles.
				
				
				{\bf Li, Okamura,  and Dohi \cite{li2019parameter}:}
				suppose an Mt/M/1/K queueing system whose job arrival follows a Non-homogeneous
				Poisson Process (NHPP) and propose a parameter estimation method for the NHPP approximately from the
				utilization data based on the maximum likelihood estimation (MLE) via the expectation maximization (EM)
				algorithm.
				
				{\bf Mei,  Gu,  Chung,   Li,  and Tang \cite{mei2019bayesian}:} A Bayesian approach for estimating vehicle queue lengths at signalized intersections using probe vehicle data.
				
				{\bf Ravner, Boxma, and  Mandjes \cite{2019Ravner}:} Develops estimation schemes for a Lévy-driven queue by sampling the workload process at Poisson times.
				
				{\bf Tan,  Yao,  Tang, and Sun  \cite{tan2019cycle}:} Cycle-based queue length estimation by fusing real-time and historical probe vehicle trajectory data, through a statistical parameter estimation method, i.e., maximum likelihood estimation (MLE)
				
				
				{\bf Zhao,  Zheng,  Wong,  Wang,  Meng, and Liu \cite{zhao2019various}:} Various methods for queue length and traffic volume estimation using probe vehicle trajectories. \\

				\medskip
				{\bf  {\color{mycolor2}2020} }
				\medskip
				
				
				{\bf Cruz,  Almeida, D’Angelo, and van Woensel \cite{cruz2018traffic}:}
				A bias-corrected version of MLE estimator of traffic intensity by the nonparametric bootstrap method
				for small and moderate samples.
				
				{\bf Tan,  Liu, Wu,  Cao, and Tang  \cite{tan2020fuzing}:} Applying probe vehicle trajectory, Bayesian theory, and Fuzing license plate recognition data and vehicle trajectory data for lane-based queue length estimation at signalized intersections.
				
				{\bf Van Phu and Farhi \cite{van2020estimation}:} Estimation of urban traffic state with probe vehicles.

				{\bf Zhang,  Liu, Chen, Yu, and Wang \cite{8759936}:} Proposes a cycle-based end-of-queue estimation method for queue length using sampled vehicle trajectory data under relatively low probe vehicle penetration rates. 
				
				{\bf Wang,  Huang,  and Lo \cite{wang2020combining}:}
				This paper from the field of traffic engineering Combines shockwave analysis and Bayesian networks for traffic parameter estimation at signalized intersections by considering queue spillback.\\
				
				\medskip
				{\bf  {\color{mycolor2}2021} }
				\medskip
				
				
				{\bf Asanjarani,  Nazarathy,  and Taylor \cite{asanjarani2021survey}:} A broad literature survey of parameter and state estimation for queueing systems.
				
				{\bf Bassamboo and Ibrahim \cite{bassamboo2021general}:} 
				Using a combination
				of queueing-theoretic analysis, real-life data analysis, and simulation, the performance of static
				and dynamic announcements are analysed, and an appropriate weighted average of them is derived.
				
				{\bf Basak and Choudhury \cite{basak2021bayesian}:}
				Finding a Bayes estimator of traffic intensity for an M/M/1 queueing model using data on queue size (number of customers present in the queue) observed at any random point in time.
				
				{\bf Comert, Amdeberhan,  Begashaw,  and Chowdhury \cite{comert2021combinatorial}:} A Combinatorial Approach for Nonparametric Short-Term Estimation of Queue
				Lengths using Probe Vehicles.

				{\bf Cruz,  Santos,   Oliveira,  and Quinino \cite{cruz2021estimation}:}
				Estimation of performance measures in a general bulk-arrival Markovian multi-server finite queue.

				{\bf Dieleman \cite{dieleman2021data}:}
				The method of MLE is used in combination with Stochastic Approximation (SA) to calibrate the arrival parameter of a G/G/1 queue via waiting time data. 
				
				{\bf Ebert, Dutta, Mengersen, Mira, Ruggeri, and Wu \cite{ebert2021likelihood}:}
				Presents a likelihood-free parameter estimation for dynamic queueing networks with a case study of passenger flow in an international airport terminal.
				
				\textbf{Krishnasamy, Sen,   Johari,  and Shakkottai \cite{krishnasamy2021learning}:} Considers regret analysis of a server allocation problem where service rates of servers are unknown. This analysis is in the context of multi-armed bandits. The work extends previous work: \cite{krishnasamy2016regret}.
				
				{\bf Lin, He, and Pang  \cite{lin2021queuing}:} Queuing network topology inference using passive end-to-end measurements originated by a single source.
				
				{\bf Mandjes and Ravner  \cite{mandjes2021hypothesis}:} In this paper the authors devise hypothesis testing procedures for L{\'e}vy-driven storage systems by sampling of the storage level.

				{\bf Singh, Acharya, Cruz,  and Quinino \cite{singh2021bayesian}:}
				Proposed a methodology to determine the sample size for an  queueing system under the Bayesian setup by observing the number of customer arrivals during the service time of a customer.
				
				{\bf Walton and Xu \cite{walton2021learning}:} Provides an extensive review of reinforcement learning ideas with a view of queueing network control. The paper also connects ideas of adversarial learning to queuing network concepts framed in the context of information uncertainty. 
				
				{\bf  Wang and Honnappa \cite{wang2021calibrating}:} Studies inference for a Cox $/~G~/~\infty$ queue, sampled at discrete time points. Uses approximate inference for maximizing a lower bound for the associated finite dimensional distribution.
				
				{\bf Zhao, Shen, and Liu \cite{zhao2021hidden}:}
				A hidden Markov model for the estimation of correlated queues in probe vehicle environments.
				
				{\bf Zhao, Wong, Zheng, and Liu \cite{zhao2021maximum}:} Maximum likelihood estimation of probe vehicle penetration rates and queue length distributions from probe vehicle data.\\ 
				
				\medskip
				{\bf  {\color{mycolor2}2022} }
				\medskip
				
				{\bf Antunes, Jacinto, and Pacheco \cite{antunes2022statistical}:}
				A discussion on statistical inference in queueing networks with probing information.
				
				
				{\bf Comert and Bagashaw \cite{COMERT2022283}:} Cycle-to-cycle queue length estimation from connected vehicles with filtering on primary parameters. This study   improve accuracy of by enhancing the low level parameter estimators using filtering algorithms.
				
				{\bf Li, Zheng,  Okamura,  and Dohi \cite{li2022hierarchical}:}
				Hierarchical Bayesian Parameter Estimation of Queueing Systems using Utilization Data.
				
				{\bf Li, Zheng,  Okamura,  and Dohi \cite{li2022parameter}:}
				Parameter estimation of a MAP/M/1/K queueing system using utilization data. In particular, the parameters are estimated by using the maximum likelihood estimation (MLE) method. 
				
				{\bf Ravner \cite{ravner2022queue}:}
				Considers an $M/G/1$ queue for which the workload process is observed periodically. The goal is to estimate the arrival rate $\lambda$
				and the parameters of the job-size distribution $G$.
				
				{\bf Singh,  Acharya,  Cruz,  and Quinino \cite{singh2022bayesian}:}
				Bayesian inference and prediction in a  queueing system.
				
				{\bf Zhong, Birge,  and Ward \cite{zhong2022learning}:} Considers a scheduling policy in time-varying multiclass many server queues with abandonment and proposes a Learn-Then-Schedule algorithm composed of a learning phase and an exploitation phase.\\
				
				\medskip
				{\bf  {\color{mycolor2}2023} }
				\medskip
				
				{\bf Chen, Liu, and Hong \cite{chen2023online}:} Proposes and studies an on-line learning algorithm for optimal control and pricing of a GI/GI/1 queue. 
				
				{\bf  Bura and Sharma \cite{bura2023maximum}:}
				Maximum likelihood and Bayesian estimation for an M/M/1 queueing model with balking.
				
				{\bf  Carmeli, Yom-Tov and Boxma \cite{carmeli2023state}:} Analyzes 
				fork-join networks and incorporates data driven estimation for such models. 
				
				
				{\bf  Inoue, Ravner, and Mandjes \cite{inoue2023estimating}:}
				Deals with parameter estimation of a queueing system with impatient customers and balking. A novel algorithm for the estimation of customer impatience is proposed and analyzed.
				
				{\bf  Luo, Deng, Chen, et. al. \cite{luo2023queue}:}
				This paper from the field of traffic engineering focuses on queue length estimation based on probe vehicle data at signalized intersections using a shockwave approach in the model.
				
				{\bf  Ravner and Wang \cite{ravner2023estimating}:}
				Parameter estimation of a queueing system where customers behave strategically, are able to choose their arrival times, and are at a Nash equilibrium. \\
				

				
				
			\newpage
				\section{Literature Analysis (up to 2017)}
				\label{sec:overview}
				
				This section was last updated in 2017 and captures a classification of the papers up to 2017 via several categories. It also serves as a background companion to our survey paper,~\cite{asanjarani2021survey}.

				
				Up to 2017, the majority of the references were in the format of journal and/or conference articles. A few are surveys, textbooks, book chapters, Ph.D. theses, and significant related materials which are listed in the table below.
				
				\begin{table}[h!]
					\begin{center}
						\begin{tabular}{ |>{\centering\arraybackslash}  m{6cm} |>{\centering\arraybackslash}  m{11cm} | }
							\hline
							\textbf{Surveys}
							&  
							\begin{center}
								\cite{bhat1987statistical}
								\cite{bhat1997statistical} and the more recent \cite{asanjarani2021survey}
							\end{center}
							\\
							\hline
							\textbf{Textbooks} 
							&
							\begin{center}
								\cite{borovkov2012stochastic}
								\cite{meyn2012markov}
								\cite{nelson2012stochastic}
								\cite{prabhu2012stochastic}
								\cite{bhat2015introduction}
								\cite{yu2015hidden}
							\end{center}
							\\
							\hline
							\textbf{Book chapters} 
							&
							\begin{center}
								\cite{bhat2019maximum}~--Chapter~2,
								\cite{bhat2008introduction}~--Chapter~10, 
								\cite{gross1998fundamentals}~--Section~6.7,
								\cite{fu2012conditional}~--Chapter~6,
								\cite{walrand1988introduction}~--Chapter~10
								
							\end{center}
							\\
							\hline
							\textbf{Ph.D. theses} 
							&
							\begin{center}
								\cite{greenberg1964parameter}
								\cite{ramirez2008bayesian-thesis}
								\cite{ibrahim2010realtimedelay}
								\cite{kauffmann2011inverse}
								\cite{horng2013inferring}
								\cite{zhao2020traffic}
							\end{center}
							\\
							\hline
							\textbf{Significant related materials} 
							&
							\begin{center}
								\cite{billingsley1961statistical} (A book)
							\end{center}
							\\
							\hline 
						\end{tabular}
					\end{center}
				\end{table}
				\subsection{Classification by Model}
				\begin{table}[h]
					\begin{center}
						\begin{tabular}{ |>{\centering\arraybackslash}  m{6cm} |>{\centering\arraybackslash}  m{11cm} | }
							\hline
							\textbf{ Model}& \textbf{References} \\
							\hline
							\textbf{M/M/1}
							& 
							\begin{center}
								\cite{clarke1957maximum}
								\cite{greenberg1967behavior}
								\cite{daley1968serial}
								\cite{jenkins1972relative}
								\cite{aigner1974parameter}
								\cite{schruben1982some}
								\cite{mcgrath1987subjective} 
								\cite{mcgrath1987subjectiveII}
								\cite{kumar1992average}
								\cite{armero1994prior}\cite{armero1994bayesian}
								\cite{masuda1995exploiting}
								\cite{rodrigues1998note}
								\cite{sohn2002robust}
								\cite{choudhury2008bayesian}
								\cite{sutton2008probabilistic}
								\cite{dey2008note}
								\cite{chen2010simulation}
								\cite{xu2010hypothesis}
								\cite{manoharan2011markovian}
								\cite{srinivas2011estimation}
								\cite{ren2012bayes}
								\cite{ren2012bayesB}
								\cite{chen2015cusum}
								\cite{quinino2016bayesian}
								\cite{zheng2000some}
							\end{center}
							\\
							\hline
							\textbf{M/M/2} {\footnotesize (Heterogenous Servers)}
							&
							\begin{center}
								\cite{dave1980maximum}
								\cite{edelman1984comments}
							\end{center} 
							\\
							\hline    
							\textbf{M/M/1/K}
							&
							\begin{center}
								\cite{mcgrath1987subjective}  \cite{mcgrath1987subjectiveII}
								\cite{alouf2001inferring}
							\end{center}   
							\\   
							\hline
						\end{tabular}
					\end{center}
				\end{table}
				\newpage
				\begin{table}[h!]
					\begin{center}
						\begin{tabular}{ |>{\centering\arraybackslash}  m{6cm} |>{\centering\arraybackslash}  m{11cm} | }
							\hline
							\textbf{ Model}& \textbf{References} \\
							\hline
							\textbf{M/M/1/$\infty$ (FIFO)}
							&
							\begin{center}
								\cite{armero1985bayesian}
							\end{center}
							\\
							\hline 
							\textbf{M/M/c}
							&
							\begin{center}
								\cite{machihara1984carried}
								\cite{hantler1989optimal}
								\cite{mcvinishconstructing}
								\cite{wang2006maximum}
								\cite{ross2007estimation}
								\cite{mcvinishconstructing}
								\cite{acharya2013maximum}
								\cite{yom2014erlang}
								\cite{chen2015cusum}
							\end{center}  
							\\
							\hline 
							\textbf{M/M/c/N}
							&
							\begin{center}
								\cite{estimates2004kuo}
								\cite{morales2007bayesian}
							\end{center}  
							\\
							\hline 
							
							\textbf{M/M/$\infty$}
							&
							\begin{center}
								\cite{armero1997bayesian}
								\cite{mandjes2009queueing}
							\end{center} 
							\\
							\hline
							\textbf{k-Par/M/1} \begin{footnotesize}
								(k-Par denotes a mixture
								of k Pareto distributions)
								
							\end{footnotesize}&
							\begin{center}
								\cite{ramirez2008bayesian}
							\end{center} 
							\\
							\hline
							\textbf{M/D/1}
							&
							\begin{center}
								\cite{novak2009determining}
								\cite{srinivas2015ml}
							\end{center}
							\\
							\hline
							\textbf{M/E$_k$/1}
							&
							\begin{center}
								\cite{jain1991confidence}
								\cite{zhang2016generalized} 
							\end{center}   
							
							\\
							\hline 
							\textbf{M/G/$\infty$} 
							{\footnotesize  (Random translation models in general)}
							&
							\begin{center}
								\cite{benes1957sufficient}
								\cite{brown1970m}
								\cite{brillinger1974cross}
								\cite{ramalhoto1987some}
								\cite{pickands1997estimation}
								\cite{bingham1999non}
								\cite{park2007choice}
								\cite{mandjes2009queueing}
								\cite{grubel2011matchmaking}
								\cite{mandjes2011m}
								\cite{nelgabats2012}
								\cite{weerasinghe2013abandonment}
								\cite{keilson1994networks}
							\end{center}
							\\
							\hline 
							\textbf{M/G/1}
							&
							\begin{center}
								\cite{daley1968serial}
								\cite{thiagarajan1979statistical}
								\cite{hernandez1983adaptive}
								\cite{hernandez1984optimal}
								\cite{jain1992relative}
								\cite{harishchandra1988note}
								\cite{gaver1990inference}  
								\cite{manjunath1996passive}
								\cite{den2017convergence}
								\cite{dimitrijevic1996inferring}
								\cite{insua1998bayesian}
								\cite{bingham1999nonparametric}
								\cite{jain2000autoregressive}
								\cite{ausin2004bayesian}
								\cite{fearnhead2004filtering}
								\cite{castellanos2006bayesian}
								\cite{rodrigo2006estimators}
								\cite{chu2007interval}
								\cite{ke2008analysis}
								\cite{ke2008estimation}
								\cite{kiessler2009technical}
								\cite{nam2009estimation}
								\cite{zhang2010confidence}
								\cite{chen2011performance}
								\cite{mohammadi2012bayesian}
								\cite{antunes2014probing}
								\cite{chen2015cusum}
								\cite{goldenshluger2015nonparametric}
								
							\end{center} 
							\\
							\hline
							\textbf{M/GI/1}
							&
							\begin{center}
								\cite{abramovstatistical}
							\end{center}  
							\\
							\hline
							\textbf{E$_k$/G/1}
							&
							\begin{center}
								\cite{fearnhead2004filtering}
							\end{center} 
							
							\\
							\hline
							
							\textbf{M/G/c}
							&
							\begin{center}
								\cite{daley1993two}
							\end{center} 
							\\
							\hline
						\end{tabular}
					\end{center}
				\end{table}
				\begin{table}[h!]
					\begin{center}
						\begin{tabular}{ |>{\centering\arraybackslash}  m{6cm} |>{\centering\arraybackslash}  m{11cm} | }
							\hline
							\textbf{ Model}& \textbf{References} \\
							\hline
							
							\textbf{M/G/c/C}
							&
							\begin{center}
								\cite{huang2001estimation}
							\end{center} 
							
							\\
							\hline   
							\textbf{G/G/1}
							&
							\begin{center}
								\cite{toyoizumi1997sengupta} 
								\cite{rodrigo1999large} 
								\cite{chu2007mean}
								\cite{chen2015cusum}
							\end{center} 
							\\
							\hline   
							
							\textbf{GI/M/1} {\footnotesize (state-dependent arrival rate)}
							&
							\begin{center}
								\cite{daley1968monte}
								\cite{pakes1971serial}
								\cite{jain1988statistical}
							\end{center}  
							\\
							\hline
							\textbf{GI/M/1/n}
							&
							\begin{center}
								\cite{abramovstatistical}
							\end{center}  
							\\
							\hline
							\textbf{GI/M/c}  
							&
							\begin{center}
								\cite{woodside1984optimal}
								\cite{ausin2007bayesian}
								\cite{ibrahim2009real}
							\end{center}
							\\
							\hline  
							\textbf{GI/G/1} 
							&
							\begin{center}
								\cite{daley1968serial}
								\cite{basawa1988large}
								\cite{pitts1994nonparametric} 
								\cite{basawa1996maximum}
								\cite{acharya1999normal}
								\cite{jain2000statistical}
								\cite{conti2002queueing}
								\cite{basawa2008parameter}
							\end{center} 
							\\
							\hline
							\textbf{GI/G/2} 
							&
							\begin{center}
								\cite{jang2001new}
							\end{center}  
							\\
							\hline
							\textbf{GI/G/c}
							&
							\begin{center}
								\cite{kim2008new}
							\end{center}
							\\
							\hline
							\textbf{G/G/$\infty$}
							&  
							\begin{center}
								\cite{ramalhoto1987some}
							\end{center} 
							\\
							\hline  
							\textbf{GI/K/1} {\footnotesize (service time has Matrix-Exponential distribution)}
							&
							\begin{center}
								\cite{asmussen1992phase}
								\cite{asmussen1992renewal}
								\cite{sengupta1989markov}
								\cite{smith1953distribution}
							\end{center} 
							\\
							\hline
							\textbf{GI/GI/s}
							&
							\begin{center}
								\cite{dong2015stochastic}
							\end{center}
							\\
							\hline   
							\textbf{GI/GI/$\infty$}
							&
							\begin{center}
								\cite{dong2015stochastic}
							\end{center}  
							\\
							\hline 
							\textbf{M/D/$\infty$} {\footnotesize (Deterministic service)}
							&
							\begin{center}
								\cite{grassmann1981technical}
							\end{center}  
							\\
							\hline
							\textbf{M/D/1}
							&  
							\begin{center}
								\cite{machihara1984carried}
							\end{center}
							\\
							\hline
							\textbf{G/D/1} 
							&
							\begin{center}
								\cite{heckmuller2009reconstructing}
								\cite{heckmueller2010reconstructing}
								
							\end{center}
							\\
							\hline
						\end{tabular}
					\end{center}
				\end{table}
				\newpage
				\begin{table}[h!]
					\begin{center}
						\begin{tabular}{ |>{\centering\arraybackslash}  m{6cm} |>{\centering\arraybackslash}  m{11cm} | }
							\hline
							\textbf{ Model}& \textbf{References} \\
							\hline
							\textbf{MAP/D/1} 
							&
							\begin{center}
								\cite{casale2008interarrival}
							\end{center}
							\\
							\hline 
							\textbf{Erlang Loss System}
							&
							\begin{center}
								\cite{ross2005estimating}
							\end{center}
							\\
							\hline
							\textbf{Jackson Networks or More General Networks
							}&
							\begin{center}
								\cite{thiruvaiyaru1991estimation}
								\cite{baccelli2009inverse}
							\end{center}
							
							\\
							\hline
							\textbf{Systems Where Little's Law Holds}
							&
							\begin{center}
								\cite{glynn1989indirect}
							\end{center}
							\\
							\hline
						\end{tabular}
					\end{center}
				\end{table}
				

				\subsection{Classification by Sampling Regime}
				\begin{table}[h!]
					\begin{center}
						\begin{tabular}{ |>{\centering\arraybackslash}  m{6cm} |>{\centering\arraybackslash}  m{9cm} | }
							\hline
							\textbf{ Sampling Regime}& \textbf{References} \\
							\hline
							\textbf{Full Observation}
							& 
							\begin{center}
								\cite{prieger2005estimation}
							\end{center}
							\\
							\hline
							\textbf{Observation at Discrete Points}
							&
							\begin{center}
								\cite{den2017convergence}
								\cite{mcvinishconstructing}
								\cite{ross2007estimation}  
							\end{center} 
							\\
							\hline    
							\textbf{Probing}
							&
							\begin{center}
								\cite{sharma1998estimating}
								\cite{alouf2001inferring}
								\cite{antunes2016estimation}
								\cite{baccelli2009inverse}
								\cite{chen1994parameter}
								\cite{comert2009queue}
								\cite{heckmuller2009reconstructing}
								\cite{hei2006model}
								\cite{hei2005light}
								\cite{kauffmann2012inverse}
								\cite{kim2017data}
								\cite{nam2009estimation}
								\cite{novak2009determining}
								\cite{antunes2014probing}
								
							\end{center}   
							\\
							\hline 
							\textbf{Queue Inference Engine}
							&
							\begin{center}
								\cite{larson1990queue}
								\cite{gawlick1990estimating}
								\cite{chandrs1994transactional}
								\cite{dimitrijevic1996inferring}
								\cite{mandelbaum1998estimating}
								\cite{frey2010queue}
								\cite{park2011analysis}
								\cite{daley1997estimating}
								\cite{daley1998moment}
								\cite{heckmuller2011reconstructing}
							\end{center}  
							\\
							\hline 
						\end{tabular}
					\end{center}
				\end{table}
				
				\newpage
				
				\subsection{Classification by Statistical Paradigm}
				\begin{table}[h!]
					\begin{center}
						\begin{tabular}{ |>{\centering\arraybackslash}  m{6cm} |>{\centering\arraybackslash}  m{9cm} | }
							\hline
							\textbf{ Statistical Paradigm}& \textbf{References} \\
							\hline
							\textbf{Bayesian}
							& 
							\begin{center}
								\cite{muddapur1972bayesian}
								\cite{reynolds1973estimating}
								\cite{warfield1984application}
								\cite{mcgrath1987subjective}
								\cite{mcgrath1987subjectiveII}
								\cite{chandrs1994transactional}
								\cite{thiruvaiyaru1992empirical} \cite{singpurwalla1992discussion}
								\cite{armero1994bayesian}
								\cite{armero1994prior}
								\cite{armero1994bayesianInf}
								\cite{armero1997bayesian}
								\cite{insua1998bayesian}
								\cite{rodrigues1998note}
								\cite{conti1999large}
								\cite{jones1999inferring}
								\cite{armero2000prediction}
								\cite{sohn2002robust}
								\cite{ausin2004bayesian}
								\cite{conti2004bootstrap}
								\cite{castellanos2006bayesian} 
								\cite{ausin2007bayesian}
								\cite{morales2007bayesian}
								\cite{choudhury2008bayesian}
								\cite{ausin2008bayesian}
								\cite{ke2009comparison}
								\cite{ramirez2010bayesian}
								\cite{sutton2011bayesian}
								\cite{mohammadi2012bayesian}
								\cite{ren2012bayes}
								\cite{ren2012bayesB}
								\cite{quinino2016bayesian}
								
							\end{center}
							\\
							\hline
							\textbf{Maximum Entropy}
							&
							\begin{center}
								\cite{he2014maximum}
							\end{center} 
							\\
							\hline    
							\textbf{Emphasis on the way of selecting sampling time}
							&
							\begin{center}
								\cite{basawa1988large} 
							\end{center}   
							\\
							\hline 
							\textbf{Non-parametric}
							&
							\begin{center}
								\cite{conti1999large}
								\cite{bingham1999non}
								\cite{jones1999inferring}
								\cite{conti2002nonparametric}
								\cite{conti2004bootstrap}
								\cite{ke2006nonparametric}
								\cite{park2007choice}
								\cite{mccabe2011efficient}
								\cite{goldenshluger2015nonparametric}
							\end{center}  
							\\
							\hline 
							\textbf{Change point detection}
							&
							\begin{center}
								\cite{jain1995estimating}
							\end{center}
							\\
							\hline 
							\textbf{Adaptive Control}
							&
							\begin{center}
								\cite{hernandez1983adaptive}
								\cite{weerasinghe2013abandonment}
							\end{center}
							\\
							\hline 
							\textbf{Sequential Inference}
							&
							\begin{center}
								\cite{jain1989problem}
								\cite{basawa1992sequential}
								\cite{vorobeychikov2015cusum}
							\end{center}
							\\
							\hline 
							\textbf{Perturbation analysis} 
							&
							\begin{center}
								\cite{ho1997perturbation}
							\end{center}
							\\
							\hline
							\textbf{Large Deviations} 
							&
							\begin{center}
								\cite{glynn1997parametric}
							\end{center} 
							\\
							\hline
						\end{tabular}
					\end{center}
				\end{table}
				
				\newpage
				
				\subsection{Classification by Application}
				\begin{table}[h!]
					\begin{center}
						\begin{tabular}{ |>{\centering\arraybackslash}  m{7.5cm} |>{\centering\arraybackslash}  m{9cm} | }
							\hline
							\textbf{ Statistical Paradigm}& \textbf{References} \\
							\hline
							\textbf{Telephone Call Centres}
							& 
							\begin{center}
								\cite{brown2005statistical}
								\cite{gorst2009asymptotic}
								\cite{ibrahim2009real}
								\cite{gans2010service}
								\cite{kim2014call}
								\cite{weerasinghe2013abandonment}
								\cite{azriel2014erlang}
							\end{center}
							\\
							\hline
							\textbf{Manufacturing}
							&
							\begin{center}
								\cite{jang1994waiting}
								\cite{castellanos2006bayesian}
								\cite{ibrahim2009real}
								\cite{chen2010simulation}
								\cite{weerasinghe2013abandonment}
							\end{center} 
							\\
							\hline    
							\textbf{Health Care}
							&
							\begin{center}
								\cite{kim2014call}
								\cite{kim2017data}
								\cite{mohammadi2012bayesian}
								\cite{whitt2017data}
								\cite{yom2014erlang}
							\end{center}   
							\\
							\hline 
							\textbf{Transportation}
							&
							\begin{center}
								\cite{castellanos2006bayesian}
								\cite{li2013freeway}
								\cite{liu2015improving}
								\cite{sutarto2015modeling}
								\cite{zammit2016joint}
								\cite{anusha2016model}
							\end{center} 
							\\
							\hline 
							\textbf{Economics}
							&
							\begin{center}
								\cite{prieger2005estimation} 
							\end{center} 
							\\
							\hline 
							\textbf{ATM}
							&
							\begin{center}
								\cite{duffield1995entropy}
								\cite{dimitrijevic1996inferring}
								\cite{conti2002nonparametric}
								\cite{conti2002queueing}
								\cite{castellanos2006bayesian}
							\end{center}
							\\
							\hline 
							\textbf{Communication/Telecommunication 
								Networks}
							&
							\begin{center}
								\cite{gawlick1990estimating}
								\cite{baccelli2009inverse}
								\cite{bakholdina2015optimal}
								\cite{massey1996estimating}
								\cite{mohajerzadeh2015efficient}
							\end{center}
							\\
							\hline 
							\textbf{Network Traffic Modelling}
							&
							\begin{center}
								\cite{sharma1999estimatingB}
								\cite{mandjes2005inferring}
								\cite{casale2008interarrival}
								\cite{liu2009real}
								\cite{ban2011real}
								\cite{sousa2011suitability}
								\cite{sutton2011bayesian}
								
							\end{center}
							\\
							\hline 
						\end{tabular}
					\end{center}
				\end{table}
				
				\newpage
				\bibliographystyle{plain}
				\bibliography{BooksChaptersDB,papersDB,ThesesisDB}

\begin{thebibliography}{100}

\bibitem{abramovstatistical}
V.M. Abramov.
\newblock Statistical analysis of single-server loss queueing systems.
\newblock {\em Methodology and Computing in Applied Probability}, pages 1--19.

\bibitem{acharya2013maximum}
S.~K. Acharya, S.V. Rodr{\'\i}guez-S{\'a}nchez, and C.E.
  Villarreal-Rodr{\'\i}guez.
\newblock Maximum likelihood estimates in an {M/M/c} queue with heterogeneous
  servers.
\newblock {\em International Journal of Mathematics in Operational Research},
  5(4):537--549, 2013.

\bibitem{acharya1999normal}
S.K. Acharya.
\newblock On normal approximation for maximum likelihood estimation from single
  server queues.
\newblock {\em Queueing Systems}, 31(3):207--216, 1999.

\bibitem{aigner1974parameter}
D.J. Aigner.
\newblock Parameter estimation from cross-sectional observations on an
  elementary queuing system.
\newblock {\em Operations Research}, 22(2):422--428, 1974.

\bibitem{almeida2018note}
M.A.C. Almeida and F.R.B. Cruz.
\newblock A note on bayesian estimation of traffic intensity in single-server
  {Markovian} queues.
\newblock {\em Communications in Statistics-Simulation and Computation},
  47(9):2577--2586, 2018.

\bibitem{alouf2001inferring}
S.~Alouf, P.~Nain, and D.~Towsley.
\newblock Inferring network characteristics via moment-based estimators.
\newblock In {\em INFOCOM 2001. Twentieth Annual Joint Conference of the IEEE
  Computer and Communications Societies. Proceedings. IEEE}, volume~2, pages
  1045--1054. IEEE, 2001.

\bibitem{amani2011multi}
P.~Amani, M.~Kihl, and A.~Robertsson.
\newblock Multi-step ahead response time prediction for single server queuing
  systems.
\newblock In {\em Computers and Communications (ISCC), 2011 IEEE Symposium on},
  pages 950--955. IEEE, 2011.

\bibitem{amini2016queue}
Z.~Amini, R.~Pedarsani, A.~Skabardonis, and P.~Varaiya.
\newblock Queue-length estimation using real-time traffic data.
\newblock In {\em Intelligent Transportation Systems (ITSC), 2016 IEEE 19th
  International Conference on}, pages 1476--1481. IEEE, 2016.

\bibitem{an2018real}
C.~An, Y.J. Wu, and W.~Xia, J.and~Huang.
\newblock Real-time queue length estimation using event-based advance detector
  data.
\newblock {\em Journal of Intelligent Transportation Systems}, 22(4):277--290,
  2018.

\bibitem{antunes2014probing}
N.~Antunes, G.~Jacinto, and A.~Pacheco.
\newblock Probing a {M/G/1} queue with general input and service times.
\newblock {\em ACM SIGMETRICS Performance Evaluation Review}, 41(3):34--36,
  2014.

\bibitem{antunes2022statistical}
N.~Antunes, G.~Jacinto, and A.~Pacheco.
\newblock Statistical inference in queueing networks with probing information.
\newblock {\em Queueing Systems}, 100(3-4):493--495, 2022.

\bibitem{antunes2016estimation}
N.~Antunes, G.~Jacinto, A.~Pacheco, and C.~Wichelhaus.
\newblock Estimation of the traffic intensity in a piecewise-stationary mt/gt/1
  queue with probing.
\newblock {\em ACM SIGMETRICS Performance Evaluation Review}, 44(2):3--5, 2016.

\bibitem{anusha2016model}
S.P. Anusha, A.~Sharma, L.~Vanajakshi, S.C. Subramanian, and L.R. Rilett.
\newblock Model-based approach for queue and delay estimation at signalized
  intersections with erroneous automated data.
\newblock {\em Journal of Transportation Engineering}, 142(5):04016013, 2016.

\bibitem{armero1985bayesian}
C.~Armero.
\newblock {Bayesian analysis of M/M/1/$\infty$/FIFO queues}.
\newblock {\em Bayesian statistics}, 2:613--618, 1985.

\bibitem{armero1994bayesianInf}
C.~Armero.
\newblock Bayesian inference in {M}arkovian queues.
\newblock {\em Queueing Systems}, 15(1):419--426, 1994.

\bibitem{armero1994bayesian}
C.~Armero and M.J. Bayarri.
\newblock Bayesian prediction in {M}/{M}/1 queues.
\newblock {\em Queueing Systems}, 15(1):401--417, 1994.

\bibitem{armero1994prior}
C.~Armero and M.J. Bayarri.
\newblock Prior assessments for prediction in queues.
\newblock {\em The Statistician}, pages 139--153, 1994.

\bibitem{armero1997bayesian}
C.~Armero and M.J. Bayarri.
\newblock A bayesian analysis of a queueing system with unlimited service.
\newblock {\em Journal of Statistical Planning and Inference}, 58(2):241--261,
  1997.

\bibitem{armero2000prediction}
C.~Armero and D.~Conesa.
\newblock Prediction in {M}arkovian bulk arrival queues.
\newblock {\em Queueing Systems}, 34(1):327--350, 2000.

\bibitem{asanjarani2021survey}
A.~Asanjarani, Y.~Nazarathy, and P.~Taylor.
\newblock A survey of parameter and state estimation in queues.
\newblock {\em Queueing Systems}, 97:39--80, 2021.

\bibitem{asmussen1992phase}
S.~Asmussen.
\newblock Phase-type representations in random walk and queueing problems.
\newblock {\em The Annals of Probability}, pages 772--789, 1992.

\bibitem{asmussen1992renewal}
S.~Asmussen and M.~Bladt.
\newblock {\em Renewal Theory and Queueing Algorithms for Matric-exponential
  Distributions}.
\newblock University of Aalborg, Institute for Electronic Systems, Department
  of Mathematics and Computer Science, 1992.

\bibitem{ausin2007bayesian}
M.C. Aus{\'\i}n, R.E. Lillo, and M.P. Wiper.
\newblock Bayesian control of the number of servers in a {GI}/{M}/c queueing
  system.
\newblock {\em Journal of Statistical Planning and Inference},
  137(10):3043--3057, 2007.

\bibitem{ausin2004bayesian}
M.C. Aus{\'\i}n, M.P. Wiper, and R.E. Lillo.
\newblock Bayesian estimation for the {M}/{G}/1 queue using a phase-type
  approximation.
\newblock {\em Journal of Statistical Planning and Inference},
  118(1-2):83--101, 2004.

\bibitem{ausin2008bayesian}
M.C. Aus{\'\i}n, M.P. Wiper, and R.E. Lillo.
\newblock Bayesian prediction of the transient behaviour and busy period in
  short-and long-tailed {GI}/{G}/1 queueing systems.
\newblock {\em Computational Statistics \& Data Analysis}, 52(3):1615--1635,
  2008.

\bibitem{azriel2014erlang}
D.~Azriel, P.D. Feigin, and A.~Mandelbaum.
\newblock Erlang s: A data-based model of servers in queueing networks.
\newblock Technical report, Working paper, 2014.

\bibitem{baccelli2009inverse}
F.~Baccelli, B.~Kauffmann, and D.~Veitch.
\newblock Inverse problems in queueing theory and internet probing.
\newblock {\em Queueing Systems}, 63(1):59--107, 2009.

\bibitem{baccelli2009towards}
F.~Baccelli, B.~Kauffmann, and D.~Veitch.
\newblock Towards multihop available bandwidth estimation.
\newblock {\em ACM SIGMETRICS Performance Evaluation Review}, 37(2):83--84,
  2009.

\bibitem{bakholdina2015optimal}
M.A. Bakholdina and A.M. Gortsev.
\newblock Optimal estimation of the states of modulated semi-synchronous
  integrated flow of events in condition of its incomplete observability.
\newblock {\em Applied Mathematical Sciences}, 9(29):1433--1451, 2015.

\bibitem{ban2011real}
X.J. Ban, P.~Hao, and Z.~Sun.
\newblock Real time queue length estimation for signalized intersections using
  travel times from mobile sensors.
\newblock {\em Transportation Research Part C: Emerging Technologies},
  19(6):1133--1156, 2011.

\bibitem{baras1984stability}
J.S. Baras, A.J. Dorsey, and A.M. Makowski.
\newblock Stability, parameter estimation and adaptive control for
  discrete-time competing queues.
\newblock In {\em The 23rd IEEE Conference on Decision and Control}, pages
  1149--1157. IEEE, 1984.

\bibitem{basak2021bayesian}
A.~Basak and A.~Choudhury.
\newblock Bayesian inference and prediction in single server m/m/1 queuing
  model based on queue length.
\newblock {\em Communications in Statistics-Simulation and Computation},
  50(6):1576--1588, 2021.

\bibitem{basawa1992sequential}
I.V. Basawa and B.R. Bhat.
\newblock Sequential inference for single server queues.
\newblock {\em OXFORD STATISTICAL SCIENCE SERIES in book "Queueing and Related
  Models, Bhat, Basawa"}, pages 325--325, 1992.

\bibitem{basawa1996maximum}
I.V. Basawa, U.N. Bhat, and R.~Lund.
\newblock Maximum likelihood estimation for single server queues from waiting
  time data.
\newblock {\em Queueing Systems}, 24(1):155--167, 1996.

\bibitem{basawa2008parameter}
I.V. Basawa, U.N. Bhat, and J.~Zhou.
\newblock Parameter estimation using partial information with applications to
  queueing and related models.
\newblock {\em Statistics \& Probability Letters}, 78(12):1375--1383, 2008.

\bibitem{basawa1997estimating}
I.V. Basawa, R.~Lund, and U.N. Bhat.
\newblock Estimating function methods of inference for queueing parameters.
\newblock {\em Lecture Notes-Monograph Series}, pages 269--284, 1997.

\bibitem{basawa1981estimation}
I.V. Basawa and N.U. Prabhu.
\newblock Estimation in single server queues.
\newblock {\em Naval Research Logistics Quarterly}, 28(3):475--487, 1981.

\bibitem{basawa1988large}
I.V. Basawa and N.U. Prabhu.
\newblock Large sample inference from single server queues.
\newblock {\em Queueing Systems}, 3(4):289--304, 1988.

\bibitem{bassamboo2021general}
A.~Bassamboo and R.~Ibrahim.
\newblock A general framework to compare announcement accuracy: Static vs.
  les-based announcement.
\newblock {\em Management Science}, 67(7):4191--4208, 2021.

\bibitem{benes1957sufficient}
V.E. Benes.
\newblock A sufficient set of statistics for a simple telephone exchange model.
\newblock {\em Bell System Tech. J}, 36:939--964, 1957.

\bibitem{bertsimas1992deducing}
D.J. Bertsimas and L.D. Servi.
\newblock Deducing queueing from transactional data: the queue inference
  engine, revisited.
\newblock {\em Operations Research}, 40:S217--S228, 1992.

\bibitem{bhat2008introduction}
U.N. Bhat.
\newblock {\em An introduction to queueing theory: modeling and analysis in
  applications}.
\newblock Birkhauser, 2008.

\bibitem{bhat2015introduction}
U.N. Bhat.
\newblock {\em An introduction to queueing theory: modeling and analysis in
  applications}.
\newblock Birkh{\"a}user, 2015.

\bibitem{bhat2019maximum}
U.N. Bhat and I.V. Basawa.
\newblock Maximum likelihood estimation in queueing systems.
\newblock In {\em Advances on Methodological and Applied Aspects of Probability
  and Statistics}, pages 13--30. CRC Press, 2019.

\bibitem{bhat1997statistical}
U.N. Bhat, G.K. Miller, and S.S. Rao.
\newblock Statistical analysis of queueing systems.
\newblock {\em Chapter in: Frontiers in queueing}, pages 351--394, 1997.

\bibitem{bhat1987statistical}
U.N. Bhat and S.S. Rao.
\newblock Statistical analysis of queueing systems.
\newblock {\em Queueing Systems}, 1(3):217--247, 1987.

\bibitem{billingsley1961statistical}
P.~Billingsley.
\newblock {\em Statistical inference for Markov processes}, volume~2.
\newblock University of Chicago Press Chicago, 1961.

\bibitem{billingsley1961statisticalPaper}
P.~Billingsley.
\newblock Statistical methods in markov chains.
\newblock {\em The Annals of Mathematical Statistics}, pages 12--40, 1961.

\bibitem{bingham1999non}
N.H. Bingham and S.M. Pitts.
\newblock Non-parametric estimation for the {M}/{G}/$\infty$ queue.
\newblock {\em Annals of the Institute of Statistical Mathematics},
  51(1):71--97, 1999.

\bibitem{bingham1999nonparametric}
N.H. Bingham and S.M. Pitts.
\newblock Nonparametric inference from {M}/{G}/l busy periods.
\newblock {\em Stochastic Models}, 15(2):247--272, 1999.

\bibitem{bladt2005statistical}
M.~Bladt and M.~S{\o}rensen.
\newblock Statistical inference for discretely observed {M}arkov jump
  processes.
\newblock {\em Journal of the Royal Statistical Society: Series B (Statistical
  Methodology)}, 67(3):395--410, 2005.

\bibitem{borovkov2012stochastic}
A.~Borovkov.
\newblock {\em Stochastic processes in queueing theory}, volume~4.
\newblock Springer Science \& Business Media, 2012.

\bibitem{brillinger1974cross}
D.R. Brillinger.
\newblock Cross-spectral analysis of processes with stationary increments
  including the stationary {G/G/$\infty$} queue.
\newblock {\em The Annals of Probability}, pages 815--827, 1974.

\bibitem{brown2005statistical}
L.~Brown, N.~Gans, A.~Mandelbaum, A.~Sakov, H.~Shen, S.~Zeltyn, and L.~Zhao.
\newblock Statistical analysis of a telephone call center.
\newblock {\em Journal of the American Statistical Association},
  100(469):36--50, 2005.

\bibitem{brown1970m}
M.~Brown.
\newblock An {M/G/$\infty$ } estimation problem.
\newblock {\em The Annals of Mathematical Statistics}, 41(2):651--654, 1970.

\bibitem{bura2023maximum}
G.S. Bura and H.~Sharma.
\newblock Maximum likelihood and bayesian estimation on {M/M/1} queueing model
  with balking.
\newblock {\em Communications in Statistics-Theory and Methods}, pages 1--29,
  2023.

\bibitem{cahoy2015transient}
D.O. Cahoy, F.~Polito, and V.~Phoha.
\newblock Transient behavior of fractional queues and related processes.
\newblock {\em Methodology and Computing in Applied Probability},
  17(3):739--759, 2015.

\bibitem{cao2003web}
J.~Cao, M.~Andersson, C.~Nyberg, and M.~Kihl.
\newblock Web server performance modeling using an {$M/G/1/K*$ PS} queue.
\newblock In {\em 10th International Conference on Telecommunications, 2003.
  ICT 2003.}, volume~2, pages 1501--1506. IEEE, 2003.

\bibitem{carmeli2023state}
N.~Carmeli, G.B. Yom-Tov, and O.J. Boxma.
\newblock State-dependent estimation of delay distributions in fork-join
  networks.
\newblock {\em Manufacturing \& Service Operations Management},
  25(3):1081--1098, 2023.

\bibitem{casale2008robust}
G.~Casale, P.~Cremonesi, and R.~Turrin.
\newblock Robust workload estimation in queueing network performance models.
\newblock In {\em Parallel, Distributed and Network-Based Processing, 2008. PDP
  2008. 16th Euromicro Conference on}, pages 183--187. IEEE, 2008.

\bibitem{casale2008interarrival}
G.~Casale, E.Z. Zhang, and E.~Smirni.
\newblock {\em Interarrival times characterization and fitting for Markovian
  traffic analysis}.
\newblock Citeseer, 2008.

\bibitem{castellanos2006bayesian}
M.E. Castellanos, J.~Morales, A.M. Mayoral, R.~Fried, and C.~Armero.
\newblock {\em On Bayesian design in finite source queues}.
\newblock Centro de Investigaci{\'o}n Operativa, Universidad Miguel
  Hern{\'a}ndez, 2006.

\bibitem{chandrs1994transactional}
K.~Chandrs and L.K. Jones.
\newblock Transactional data inference for telecommunication models.
\newblock In {\em presentation at First Annual Technical Conference on
  Telecommunications R \& D in Massachusetts, University of Massachusetts,
  Lowell, Massachusetts}, 1994.

\bibitem{chen1988empirical}
H.~Chen, J.M. Harrison, A.~Mandelbaum, A.~Van~Ackere, and L.M. Wein.
\newblock Empirical evaluation of a queueing network model for semiconductor
  wafer fabrication.
\newblock {\em Operations Research}, pages 202--215, 1988.

\bibitem{chen2011performance}
N.~Chen, Y.~Yuan, and S.~Zhou.
\newblock Performance analysis of queue length monitoring of {M/G/1} systems.
\newblock {\em Naval Research Logistics (NRL)}, 58(8):782--794, 2011.

\bibitem{chen2010simulation}
N.~Chen and S.~Zhou.
\newblock Simulation-based estimation of cycle time using quantile regression.
\newblock {\em IIE Transactions}, 43(3):176--191, 2010.

\bibitem{chen2015cusum}
N.~Chen and S.~Zhou.
\newblock Cusum statistical monitoring of {M/M/1} queues and extensions.
\newblock {\em Technometrics}, 57(2):245--256, 2015.

\bibitem{chen1994parameter}
T.M. Chen, J.~Walrand, and D.G. Messerschmitt.
\newblock Parameter estimation for partially observed queues.
\newblock {\em Communications, IEEE Transactions on}, 42(9):2730--2739, 1994.

\bibitem{chen2023online}
X.~Chen, Y.~Liu, and G.~Hong.
\newblock An online learning approach to dynamic pricing and capacity sizing in
  service systems.
\newblock {\em Operations Research}, 2023.

\bibitem{chick2006subjective}
S.E. Chick.
\newblock Subjective probability and bayesian methodology.
\newblock {\em Chapter in: Handbooks in Operations Research and Management
  Science}, 13:225--257, 2006.

\bibitem{choudhury2008bayesian}
A.~Choudhury and A.C. Borthakur.
\newblock Bayesian inference and prediction in the single server {M}arkovian
  queue.
\newblock {\em Metrika}, 67(3):371--383, 2008.

\bibitem{chow2013observable}
J.~Chow.
\newblock On observable chaotic maps for queuing analysis.
\newblock {\em Transportation Research Record: Journal of the Transportation
  Research Board}, (2390):138--147, 2013.

\bibitem{chu2006confidence}
Y.K. Chu and J.C. Ke.
\newblock Confidence intervals of mean response time for an {M}/{G}/1 queueing
  system: Bootstrap simulation.
\newblock {\em Applied Mathematics and Computation}, 180(1):255--263, 2006.

\bibitem{chu2007interval}
Y.K. Chu and J.C. Ke.
\newblock Interval estimation of mean response time for a {G}/{M}/1 queueing
  system: empirical laplace function approach.
\newblock {\em Mathematical Methods in the Applied Sciences}, 30(6):707--715,
  2007.

\bibitem{chu2007mean}
Y.K. Chu and J.C. Ke.
\newblock Mean response time for a {G}/{G}/1 queueing system: Simulated
  computation.
\newblock {\em Applied Mathematics and Computation}, 186(1):772--779, 2007.

\bibitem{chu2009analysis}
Y.K. Chu and J.C. Ke.
\newblock Analysis of intensity for a queueing system: bootstrapping
  computation.
\newblock {\em International Journal of Services Operations and Informatics},
  4(3):195--211, 2009.

\bibitem{clarke1957maximum}
A.B. Clarke.
\newblock Maximum likelihood estimates in a simple queue.
\newblock {\em The Annals of Mathematical Statistics}, 28(4):1036--1040, 1957.

\bibitem{comert2016queue}
G.~Comert.
\newblock Queue length estimation from probe vehicles at isolated
  intersections: Estimators for primary parameters.
\newblock {\em European Journal of Operational Research}, 252(2):502--521,
  2016.

\bibitem{comert2021combinatorial}
G.~Comert, T.~Amdeberhan, N.~Begashaw, and M.~Chowdhury.
\newblock A combinatorial approach for nonparametric short-term estimation of
  queue lengths using probe vehicles.
\newblock {\em arXiv preprint arXiv:2112.04551}, 2021.

\bibitem{COMERT2022283}
G.~Comert and N.~Begashaw.
\newblock Cycle-to-cycle queue length estimation from connected vehicles with
  filtering on primary parameters.
\newblock {\em International Journal of Transportation Science and Technology},
  11(2):283--297, 2022.

\bibitem{comert2009queue}
G.~Comert and M.~Cetin.
\newblock Queue length estimation from probe vehicle location and the impacts
  of sample size.
\newblock {\em European Journal of Operational Research}, 197(1):196--202,
  2009.

\bibitem{conti1999large}
P.L. Conti.
\newblock Large sample bayesian analysis for {Geo/G/1} discrete-time queueing
  models.
\newblock {\em The Annals of Statistics}, 27(6):1785--1807, 1999.

\bibitem{conti2002nonparametric}
P.L. Conti.
\newblock Nonparametric statistical analysis of discrete-time queues, with
  applications to atm teletraffic data.
\newblock {\em Stochastic Models}, 18(4):497--527, 2002.

\bibitem{conti2004bootstrap}
P.L. Conti.
\newblock Bootstrap approximations for bayesian analysis of {Geo/G/1}
  discrete-time queueing models.
\newblock {\em Journal of statistical planning and inference}, 120(1-2):65--84,
  2004.

\bibitem{conti2002queueing}
P.L. Conti and L.~De~Giovanni.
\newblock Queueing models and statistical analysis for atm based networks.
\newblock {\em Sankhy{\=a}: The Indian Journal of Statistics, Series B}, pages
  50--75, 2002.

\bibitem{cox1955statistical}
D.R. Cox.
\newblock The statistical analysis of congestion.
\newblock {\em Journal of the Royal Statistical Society. Series A (General)},
  118(3):324--335, 1955.

\bibitem{cox1965some}
D.R. Cox.
\newblock Some problems of statistical analysis connected with congestion.
\newblock In {\em Proc. of the Symp. on Congestion Theory}, 1965.

\bibitem{cruz2018traffic}
F.R.B. Cruz, M.A.C. Almeida, M.F.S.V. D’Angelo, and T.~van Woensel.
\newblock Traffic intensity estimation in finite markovian queueing systems.
\newblock {\em Mathematical Problems in Engineering}, 2018, 2018.

\bibitem{cruz2016bayesian}
F.R.B. Cruz, R.C. Quinino, and L.L. Ho.
\newblock Bayesian estimation of traffic intensity based on queue length in a
  multi-server {$M/M/s$} queue.
\newblock {\em Communications in Statistics-Simulation and Computation},
  (just-accepted):00--00, 2016.

\bibitem{cruz2017bayesian}
F.R.B. Cruz, R.C. Quinino, and L.L. Ho.
\newblock Bayesian estimation of traffic intensity based on queue length in a
  multi-server m/m/s queue.
\newblock {\em Communications in Statistics-Simulation and Computation},
  46(9):7319--7331, 2017.

\bibitem{cruz2021estimation}
F.R.B. Cruz, M.A.C. Santos, F.L.P. Oliveira, and R.C. Quinino.
\newblock Estimation in a general bulk-arrival markovian multi-server finite
  queue.
\newblock {\em Operational Research}, 21(1):73--89, 2021.

\bibitem{daley1968monte}
D.J. Daley.
\newblock Monte carlo estimation of the mean queue size in a stationary
  {GI/M/1} queue.
\newblock {\em Operations Research}, pages 1002--1005, 1968.

\bibitem{daley1968serial}
D.J. Daley.
\newblock The serial correlation coefficients of waiting times in a stationary
  single server queue.
\newblock {\em Journal of the Australian Mathematical Society}, 8(683-699):27,
  1968.

\bibitem{daley1992exploiting}
D.J. Daley and L.D. Servi.
\newblock Exploiting {M}arkov chains to infer queue length from transactional
  data.
\newblock {\em Journal of Applied Probability}, 29(3):713--732, 1992.

\bibitem{daley1993two}
D.J. Daley and L.D. Servi.
\newblock A two-point {M}arkov chain boundary-value problem.
\newblock {\em Advances in Applied Probability}, 25(3):607--630, 1993.

\bibitem{daley1997estimating}
D.J. Daley and L.D. Servi.
\newblock Estimating waiting times from transactional data.
\newblock {\em INFORMS Journal on Computing}, 9(2):224--229, 1997.

\bibitem{daley1998moment}
D.J. Daley and L.D. Servi.
\newblock Moment estimation of customer loss rates from transactional data.
\newblock {\em Journal of Applied Mathematics and Stochastic Analysis},
  11(3):301--310, 1998.

\bibitem{dave1980maximum}
U.~Dave and Y.K. Shah.
\newblock Maximum likelihood estimates in a {M}/{M}/2 queue with heterogeneous
  servers.
\newblock {\em Journal of the Operational Research Society}, 31(5):423--426,
  1980.

\bibitem{den2017convergence}
A.V. den Boer, M.~Mandjes, et~al.
\newblock Convergence rates of {L}aplace-transform based estimators.
\newblock {\em Bernoulli}, 23(4A):2533--2557, 2017.

\bibitem{dey2008note}
S.~Dey.
\newblock A note on bayesian estimation of the traffic intensity in {M}/{M}/1
  queue and queue characteristics under quadratic loss function.
\newblock {\em Data Science Journal}, 7(0):148--154, 2008.

\bibitem{dieleman2021data}
N.A. Dieleman.
\newblock Data-driven fitting of the {G/G/1} queue.
\newblock {\em Journal of Systems Science and Systems Engineering},
  30(1):17--28, 2021.

\bibitem{dimitrijevic1996inferring}
D.D. Dimitrijevic.
\newblock Inferring most likely queue length from transactional data.
\newblock {\em Operations Research Letters}, 19(4):191--199, 1996.

\bibitem{dinh2014architecture}
T.V. Dinh, L.L.H. Andrew, and Y.~Nazarathy.
\newblock Architecture and robustness tradeoffs in speed-scaled queues with
  application to energy management.
\newblock {\em International Journal of Systems Science}, 45(8):1728--1739,
  2014.

\bibitem{dong2015stochastic}
J.~Dong and W.~Whitt.
\newblock Stochastic grey-box modeling of queueing systems: fitting
  birth-and-death processes to data.
\newblock {\em Queueing Systems}, 79(3-4):391--426, 2015.

\bibitem{dong2015using}
J.~Dong and W.~Whitt.
\newblock Using a birth-and-death process to estimate the steady-state
  distribution of a periodic queue.
\newblock {\em Naval Research Logistics (NRL)}, 62(8):664--685, 2015.

\bibitem{doucet2006optimal}
A.~Doucet, L.~Montesano, and A.~Jasra.
\newblock Optimal filtering for partially observed point processes using
  trans-dimensional sequential monte carlo.
\newblock In {\em Acoustics, Speech and Signal Processing, 2006. ICASSP 2006
  Proceedings. 2006 IEEE International Conference on}, volume~5, pages V--V.
  IEEE, 2006.

\bibitem{duffield2000large}
N.G. Duffield.
\newblock A large deviation analysis of errors in measurement based admission
  control to buffered and bufferless resources.
\newblock {\em Queueing Systems}, 34(1):131--168, 2000.

\bibitem{duffield1995entropy}
N.G. Duffield, J.T. Lewis, N.~O'Connell, R.~Russell, and F.~Toomey.
\newblock Entropy of atm traffic streams: a tool for estimating qos parameters.
\newblock {\em Selected Areas in Communications, IEEE Journal on},
  13(6):981--990, 1995.

\bibitem{duffy2010most}
K.R. Duffy and S.P. Meyn.
\newblock Most likely paths to error when estimating the mean of a reflected
  random walk.
\newblock {\em Performance Evaluation}, 67(12):1290--1303, 2010.

\bibitem{duffy2009estimating}
K.R. Duffy and S.P. Meyn.
\newblock Estimating Loynes exponent.
\newblock {\em Queueing Systems}, 68(3-4):285--293, 2011.

\bibitem{duffy2012large}
K.R. Duffy and S.P. Meyn.
\newblock Large deviation asymptotics for busy periods.
\newblock 2012.

\bibitem{ebert2021likelihood}
A.~Ebert, R.~Dutta, K.~Mengersen, A.~Mira, F.~Ruggeri, and P.~Wu.
\newblock Likelihood-free parameter estimation for dynamic queueing networks:
  Case study of passenger flow in an international airport terminal.
\newblock {\em Journal of the Royal Statistical Society. Series C: Applied
  Statistics}, 70(3):770--792, 2021.

\bibitem{edelman1984comments}
D.B. Edelman and D.E. McKellar.
\newblock Comments on maximum likelihood estimates in a {M}/{M}/2 queue with
  heterogeneous servers.
\newblock {\em Journal of the Operational Research Society}, 35(2):149--150,
  1984.

\bibitem{edelmann2014nonparametric}
D.~Edelmann and C.~Wichelhaus.
\newblock Nonparametric inference for queueing networks of {Geom $X/G/\infty$}
  queues in discrete time.
\newblock {\em Advances in Applied Probability}, 46(3):790--811, 2014.

\bibitem{efrosinin2015confidence}
D.~Efrosinin, A.~Winkler, and P.~Martin.
\newblock Confidence intervals for performance measures of {M/M/1} queue with
  constant retrial policy.
\newblock {\em Asia-Pacific Journal of Operational Research}, 32(06):1550046,
  2015.

\bibitem{emami2019neural}
A.~Emami, M.~Sarvi, and S.~Asadi~Bagloee.
\newblock A neural network algorithm for queue length estimation based on the
  concept of k-leader connected vehicles.
\newblock {\em Journal of Modern Transportation}, 27:341--354, 2019.

\bibitem{eschenbach1984statistical}
W.~Eschenbach.
\newblock Statistical inference for queueing models.
\newblock {\em Series Statistics}, 15(3):451--462, 1984.

\bibitem{fabris2012investigation}
I.N. Fabris-Rotelli, C.~Kraamwinkel, et~al.
\newblock An investigation and historical overview of the {G/M} and {M/G}
  queueing processes.
\newblock In {\em South African Statistical Journal Proceedings: Proceedings of
  the 54th Annual Conference of the South African Statistical Association:
  Congress 1}, pages 18--25. Sabinet Online, 2012.

\bibitem{fearnhead2004filtering}
P.~Fearnhead.
\newblock Filtering recursions for calculating likelihoods for queues based on
  inter-departure time data.
\newblock {\em Statistics and Computing}, 14(3):261--266, 2004.

\bibitem{fendick1989measurements}
K.W. Fendick and W.~Whitt.
\newblock Measurements and approximations to describe the offered traffic and
  predict the average workload in a single-server queue.
\newblock {\em Proceedings of the IEEE}, 77(1):171--194, 1989.

\bibitem{feng2011estimation}
H.~Feng, P.~Dube, and L.~Zhang.
\newblock On estimation problems for the {G/G/°} queue.
\newblock {\em ACM SIGMETRICS Performance Evaluation Review}, 39(3):40--42,
  2011.

\bibitem{frey2010queue}
J.C. Frey and E.H. Kaplan.
\newblock Queue inference from periodic reporting data.
\newblock {\em Operations Research Letters}, 38(5):420--426, 2010.

\bibitem{fu2012conditional}
M.C. Fu and J.Q. Hu.
\newblock {\em Conditional Monte Carlo: Gradient estimation and optimization
  applications}, volume 392.
\newblock Springer Science \& Business Media, 2012.

\bibitem{ganesh1998bayesian}
A.~Ganesh, P.~Green, N.~O'Connell, and S.~Pitts.
\newblock Bayesian network management.
\newblock {\em Queueing Systems}, 28(1):267--282, 1998.

\bibitem{gans2010service}
N.~Gans, N.~Liu, A.~Mandelbaum, H.~Shen, H.~Ye, et~al.
\newblock Service times in call centers: Agent heterogeneity and learning with
  some operational consequences.
\newblock In {\em Borrowing Strength: Theory Powering Applications--A
  Festschrift for Lawrence D. Brown}, pages 99--123. Institute of Mathematical
  Statistics, 2010.

\bibitem{gaver1990inference}
D.P. Gaver and P.A. Jacobs.
\newblock On inference concerning time-dependent queue performance: the
  {M}/{G}/1 example.
\newblock {\em Queueing Systems}, 6(1):261--275, 1990.

\bibitem{gawlick1990estimating}
R.~Gawlick.
\newblock Estimating disperse network queues: The queue inference engine.
\newblock {\em ACM SIGCOMM Computer Communication Review}, 20(5):111--118,
  1990.

\bibitem{ghorbani2016ml}
M.~Ghorbani-Mandolakani and M.R. Salehi~Rad.
\newblock Ml and bayes estimation in a two-phase tandem queue with a second
  optional service and random feedback.
\newblock {\em Communications in Statistics-Theory and Methods},
  45(9):2576--2591, 2016.

\bibitem{glynn1993estimating}
P.W. Glynn, B.~Melamed, and W.~Whitt.
\newblock Estimating customer and time averages.
\newblock {\em Operations research}, pages 400--408, 1993.

\bibitem{glynn1997parametric}
P.W. Glynn and M.~Torres.
\newblock Parametric estimation of tail probabilities for the single-server
  queue.
\newblock {\em Chapter in: Frontiers in queueing}, pages 449--462, 1997.

\bibitem{glynn1989indirect}
P.W. Glynn and W.~Whitt.
\newblock Indirect estimation via ${L}= \lambda {W}$.
\newblock {\em Operations Research}, pages 82--103, 1989.

\bibitem{glynn28estimating}
P.W. Glynn and A.J. Zeevi.
\newblock Estimating tail probabilities in queues via extremal statistics.
  analysis of communication networks: Call centres, traffic and performance,
  135--158.
\newblock {\em Fields Inst. Commun}, 28.

\bibitem{goldenshluger2015nonparametric}
A.~Goldenshluger.
\newblock Nonparametric estimation of service time distribution in the {
  M/G/$\infty$} queue and related estimation problems.
\newblock {\em arXiv preprint arXiv:1508.00076}, 2015.

\bibitem{gordon1980impact}
K.D. Gordon and Lawrence~W. D.
\newblock The impact of certain parameter estimation errors in queueing network
  models.
\newblock In {\em Proceedings of the 1980 international symposium on Computer
  performance modelling, measurement and evaluation}, pages 3--9, 1980.

\bibitem{gorst2009asymptotic}
A.~Gorst-Rasmussen and M.B. Hansen.
\newblock Asymptotic inference for waiting times and patiences in queues with
  abandonment.
\newblock {\em Communications in Statistics--Simulation and Computation},
  38(2):318--334, 2009.

\bibitem{goyal1972maximum}
T.L. Goyal and C.M. Harris.
\newblock Maximum-likelihood estimates for queues with state-dependent service.
\newblock {\em Sankhy{\=a}: The Indian Journal of Statistics, Series A},
  34(1):65--80, 1972.

\bibitem{grassmann1981technical}
W.K. Grassmann.
\newblock Technical note. the optimal estimation of the expected number in a
  {M/D/8} queueing system.
\newblock {\em Operations Research}, 29(6):1208--1211, 1981.

\bibitem{greenberg1964parameter}
I.~Greenberg.
\newblock {\em Parameter estimation in a simple queue}.
\newblock PhD thesis, Doctoral Thesis at New York University, 1964.

\bibitem{greenberg1967behavior}
I.~Greenberg.
\newblock The behavior of a simple queue at various times and epochs.
\newblock {\em SIAM Review}, 9(2):234--248, 1967.

\bibitem{gross1998fundamentals}
D.~Gross and C.~Harris.
\newblock Fundamentals of queueing theory.
\newblock 1998.

\bibitem{grubel2011matchmaking}
R.~Gr{\"u}bel and H.~Wegener.
\newblock Matchmaking and testing for exponentiality in the {M/G/$\infty$}
  queue.
\newblock {\em Journal of Applied Probability}, 48(1):131--144, 2011.

\bibitem{gu2017traffic}
Y.~Gu, Z.~Qian, and G.~Zhang.
\newblock Traffic state estimation for urban road networks using a link queue
  model.
\newblock {\em Transportation research record}, 2623(1):29--39, 2017.

\bibitem{gurvich2013excursion}
I.~Gurvich, J.~Huang, and A.~Mandelbaum.
\newblock Excursion-based universal approximations for the erlang-a queue in
  steady-state.
\newblock {\em Mathematics of Operations Research}, 39(2):325--373, 2013.

\bibitem{halfin1982linear}
S.~Halfin.
\newblock Linear estimators for a class of stationary queueing processes.
\newblock {\em Operations Research}, 30(3):515--529, 1982.

\bibitem{hall2004nonparametric}
P.~Hall and J.~Park.
\newblock Nonparametric inference about service time distribution from indirect
  measurements.
\newblock {\em Journal of the Royal Statistical Society: Series B (Statistical
  Methodology)}, 66(4):861--875, 2004.

\bibitem{hall1991using}
S.A. Hall and R.C. Larson.
\newblock Using partial queue-length information to improve the queue inference
  engine's performance.
\newblock 1991.

\bibitem{hansen2006nonparametric}
M.B. Hansen and S.M. Pitts.
\newblock Nonparametric inference from the {M}/{G}/1 workload.
\newblock {\em Bernoulli}, 12(4):737--759, 2006.

\bibitem{hantler1989optimal}
S.L. Hantler and Z.~Rosberg.
\newblock Optimal estimation for an {M}/{M}/c queue with time varying
  parameters.
\newblock {\em Stochastic Models}, 5(2):295--313, 1989.

\bibitem{harishchandra1988note}
K.~Harishchandra and S.S. Rao.
\newblock A note on statistical inference about the traffic intensity parameter
  in {M}/{$E_k$}/1 queue.
\newblock {\em Sankhy{\=a}: The Indian Journal of Statistics, Series B}, pages
  144--148, 1988.

\bibitem{harris1973some}
C.M. Harris.
\newblock Some new results in the statistical analysis of queues.
\newblock Technical report, DTIC Document, 1973.

\bibitem{he2014maximum}
D.~He, R.~Li, Q.~Huang, and P.~Lei.
\newblock Maximum entropy principle based estimation of performance
  distribution in queueing theory.
\newblock 2014.

\bibitem{heckmueller2010reconstructing}
S.~Heckmueller and B.E. Wolfinger.
\newblock Reconstructing arrival processes to discrete queueing systems by
  inverse load transformation.
\newblock {\em Simulation}, 2010.

\bibitem{heckmuller2009reconstructing}
S.~Heckmuller and B.E. Wolfinger.
\newblock Reconstructing arrival processes to {G}/{D}/1 queueing systems and
  tandem networks.
\newblock In {\em Performance Evaluation of Computer \& Telecommunication
  Systems, 2009. SPECTS 2009. International Symposium on}, volume~41, pages
  361--368. IEEE, 2009.

\bibitem{heckmuller2011reconstructing}
S.~Heckm{\"u}ller and B.E. Wolfinger.
\newblock Reconstructing arrival processes to discrete queueing systems by
  inverse load transformation.
\newblock {\em Simulation}, 87(12):1033--1047, 2011.

\bibitem{hei2005light}
X.~Hei, B.~Bensaou, and D.H.K. Tsang.
\newblock A light-weight available bandwidth inference methodology in a
  queueing analysis approach.
\newblock In {\em Communications, 2005. ICC 2005. 2005 IEEE International
  Conference on}, volume~1, pages 120--124. IEEE, 2005.

\bibitem{hei2006model}
X.~Hei, B.~Bensaou, and D.H.K. Tsang.
\newblock Model-based end-to-end available bandwidth inference using queueing
  analysis.
\newblock {\em Computer Networks}, 50(12):1916--1937, 2006.

\bibitem{hernandez1983adaptive}
O.~Hernandez-Lerma and S.I. Marcus.
\newblock Adaptive control of service in queueing systems.
\newblock {\em Systems \& Control Letters}, 3(5):283--289, 1983.

\bibitem{hernandez1984optimal}
O.~Hern{\`a}ndez-Lerma and S.I. Marcus.
\newblock Optimal adaptive control of priority assignment in queueing systems*.
\newblock {\em Systems \& Control Letters}, 4(2):65--72, 1984.

\bibitem{heyde1992some}
C.C. Heyde.
\newblock Some results on inference for stationary processes and queueing
  systems.
\newblock {\em OXFORD STATISTICAL SCIENCE SERIES in book "Queueing and Related
  Models, Bhat, Basawa"}, pages 337--337, 1992.

\bibitem{ho1997perturbation}
Y.C. Ho and C.G. Cassandras.
\newblock Perturbation analysis for control and optimization of queueing
  systems: an overview and the state of the art.
\newblock {\em Chapter in: Frontiers in Queueing"(J. Dshalalow, Ed.), CRC
  Press}, pages 395--420, 1997.

\bibitem{horng2013inferring}
T.C. Horng.
\newblock Inferring queueing network models from high-precision location
  tracking data.
\newblock 2013.

\bibitem{hu2012parameter}
Y.~Hu and Lee.
\newblock Parameter estimation for a reflected fractional brownian motion based
  on its local time.
\newblock 2012.

\bibitem{huang2001estimation}
M.L. Huang and P.~Brill.
\newblock On estimation in {M}/{G}/c/c queues.
\newblock {\em International Transactions in Operational Research},
  8(6):647--657, 2001.

\bibitem{ibrahim2010realtimedelay}
R.~Ibrahim.
\newblock {\em Real-Time Delay Prediction in Customer Service Systems}.
\newblock PhD thesis, Doctoral Thesis at Columbia University, 2010.

\bibitem{ibrahim2009real}
R.~Ibrahim and W.~Whitt.
\newblock Real-time delay estimation based on delay history.
\newblock {\em Manufacturing \& Service Operations Management}, 11(3):397--415,
  2009.

\bibitem{ibrahim2011wait}
R.~Ibrahim and W.~Whitt.
\newblock Wait-time predictors for customer service systems with time-varying
  demand and capacity.
\newblock {\em Operations Research-Baltimore}, 59(5):1106--1118, 2011.

\bibitem{inoue2023estimating}
Y.~Inoue, L.~Ravner, and M.~Mandjes.
\newblock Estimating customer impatience in a service system with unobserved
  balking.
\newblock {\em Stochastic Systems}, 13(2):181--210, 2023.

\bibitem{insua1998bayesian}
D.R. Insua, M.~Wiper, and F.~Ruggeri.
\newblock Bayesian analysis of {$M/E_r/1$} and {$M/H_k/1$} queues.
\newblock {\em Queueing Systems}, 30(3):289--308, 1998.

\bibitem{jain1991comparison}
S.~Jain.
\newblock Comparison of confidence intervals of traffic intensity for
  {M}/{E}$_k$/1 queueing systems.
\newblock {\em Statistical Papers}, 32(1):167--174, 1991.

\bibitem{jain1992relative}
S.~Jain.
\newblock Relative efficiency of a parameter for a {M}/{G}/1 queueing system
  based on reduced and full likelihood functions.
\newblock {\em Communications in Statistics-Simulation and Computation},
  21(2):597--606, 1992.

\bibitem{jain1995estimating}
S.~Jain.
\newblock Estimating changes in traffic intensity for {M}/{M}/1 queueing
  systems.
\newblock {\em Microelectronics and Reliability}, 35(11):1395--1400, 1995.

\bibitem{jain2000autoregressive}
S.~Jain.
\newblock An autoregressive process and its application to queueing model.
\newblock {\em Metron-International Journal of Statistics}, 58(1-2):131--138,
  2000.

\bibitem{jain2000statistical}
S.~Jain and T.S.S. Rao.
\newblock Statistical inference for the {GI}/{G}/1 queue using diffusion
  approximation.
\newblock {\em International Journal of Information and Management Sciences},
  11(2):1--12, 2000.

\bibitem{jain1988statistical}
S.~Jain and J.G.C. Templeton.
\newblock Statistical inference for {G}/{M}/1 queueing system.
\newblock {\em Operations Research Letters}, 7(6):309--313, 1988.

\bibitem{jain1989problem}
S.~Jain and J.G.C. Templeton.
\newblock Problem of statistical inference to control the traffic intensity.
\newblock {\em Sequential Analysis}, 8(2):135--146, 1989.

\bibitem{jain1991confidence}
S.~Jain and J.G.C. Templeton.
\newblock Confidence interval for {M}/{M}/2 queue with heterogeneous servers.
\newblock {\em Operations Research Letters}, 10(2):99--101, 1991.

\bibitem{jang1994waiting}
J.~Jang and C.R. Liu.
\newblock Waiting time estimation in a manufacturing system using the number of
  machine idle periods.
\newblock {\em European Journal of Operational Research}, 78(3):426--440, 1994.

\bibitem{jang2001new}
J.~Jang, J.~Suh, and C.R. Liu.
\newblock A new procedure to estimate waiting time in {GI}/{G}/2 system by
  server observation.
\newblock {\em Computers \& Operations Research}, 28(6):597--611, 2001.

\bibitem{jenkins1972relative}
J.H. Jenkins.
\newblock The relative efficiency of direct and maximum likelihood estimates of
  mean waiting time in the simple queue, {M/M/1}.
\newblock {\em Journal of Applied Probability}, pages 396--403, 1972.

\bibitem{jones1999inferring}
L.K. Jones.
\newblock Inferring balking behavior from transactional data.
\newblock {\em Operations Research}, pages 778--784, 1999.

\bibitem{jones2012remarks}
L.K. Jones.
\newblock Remarks on queue inference from departure data alone and the
  importance of the queue inference engine.
\newblock {\em Operations Research Letters}, 2012.

\bibitem{jones1994efficient}
L.K. Jones and R.C. Larson.
\newblock Efficient computation of probabilities of events described by order
  statistics and applications to queue inference.
\newblock 1994.

\bibitem{kannan2014parameter}
K.S. Kannan and A.~Jabarali.
\newblock Parameter estimation of single server queue with working vacations.
\newblock {\em Research \& Reviews: Journal of Statistics (Special Issue on
  Recent Statistical Methodologies and Applications)}, 2:94--98, 2014.

\bibitem{kauffmann2011inverse}
B.~Kauffmann.
\newblock {\em Inverse problems in networks}.
\newblock PhD thesis, Universit{\'e} Pierre et Marie Curie-Paris VI, 2011.

\bibitem{kauffmann2012inverse}
B.~Kauffmann.
\newblock Inverse problems in bandwidth sharing networks.
\newblock In {\em Proceedings of the 24th International Teletraffic Congress},
  page~6. International Teletraffic Congress, 2012.

\bibitem{ke2006nonparametric}
J.C. Ke and Y.K. Chu.
\newblock Nonparametric and simulated analysis of intensity for a queueing
  system.
\newblock {\em Applied Mathematics and Computation}, 183(2):1280--1291, 2006.

\bibitem{ke2009comparison}
J.C. Ke and Y.K. Chu.
\newblock Comparison on five estimation approaches of intensity for a queueing
  system with short run.
\newblock {\em Computational Statistics}, 24(4):567--582, 2009.

\bibitem{ke2008analysis}
J.C. Ke, M.Y. Ko, and K.C. Chiou.
\newblock Analsis of 2 4 factorial design for a controllable {M}/{G}/1 system.
\newblock {\em Mathematical and Computational Applications}, 13(3):165--174,
  2008.

\bibitem{ke2008estimation}
J.C. Ke, M.Y. Ko, and S.H. Sheu.
\newblock Estimation comparison on busy period for a controllable {M}/{G}/1
  system with bicriterion policy.
\newblock {\em Simulation Modelling Practice and Theory}, 16(6):645--655, 2008.

\bibitem{keiding1975maximum}
N.~Keiding.
\newblock Maximum likelihood estimation in the birth-and-death process.
\newblock {\em The Annals of Statistics}, 3(2):363--372, 1975.

\bibitem{keilson1994networks}
J.~Keilson and L.D. Servi.
\newblock Networks of non-homogeneous {M/G/}$\infty$ systems.
\newblock {\em Journal of Applied Probability}, 31(A):157--168, 1994.

\bibitem{estimates2004kuo}
S.C.~Chen K.H.~Wang and J.C. Ke.
\newblock Maximum likelihood estimates and confidence intervals of an
  {M}/{M}/{R}/{N} queue with balking and heterogeneous servers.
\newblock {\em RAIRO -- Operations Research}, 38(3):227--241, 2004.

\bibitem{kiessler2009technical}
P.C. Kiessler and R.~Lund.
\newblock Technical note: Traffic intensity estimation.
\newblock {\em Naval Research Logistics (NRL)}, 56(4):385--387, 2009.

\bibitem{kim2012statisticalSupplementary}
S.~Kim and W.~Whitt.
\newblock Statistical analysis with little's law, supplementary material: More
  on the call center data.
\newblock {\em Supplementary Material to Preprint}, 2012.

\bibitem{kim2012statistical}
S.~Kim and W.~Whitt.
\newblock Statistical analysis with littleÕs law.
\newblock {\em preparation. Columbia University}, 2012.

\bibitem{kim2012estimatingAppendix}
S.H. Kim and W.~Whitt.
\newblock Appendix to estimating waiting times with the time-varying littleÕs
  law.
\newblock 2012.

\bibitem{kim2012estimating}
S.H. Kim and W.~Whitt.
\newblock Estimating waiting times with the time-varying littleÕs law.
\newblock 2012.

\bibitem{kim2014call}
S.H. Kim and W.~Whitt.
\newblock Are call center and hospital arrivals well modeled by nonhomogeneous
  poisson processes?
\newblock {\em Manufacturing \& Service Operations Management}, 16(3):464--480,
  2014.

\bibitem{kim2014choosing}
S.H. Kim and W.~Whitt.
\newblock Choosing arrival process models for service systems: Tests of a
  nonhomogeneous poisson process.
\newblock {\em Naval Research Logistics (NRL)}, 61(1):66--90, 2014.

\bibitem{kim2015power}
S.H. Kim and W.~Whitt.
\newblock The power of alternative kolmogorov-smirnov tests based on
  transformations of the data.
\newblock {\em ACM Transactions on Modeling and Computer Simulation (TOMACS)},
  25(4):24, 2015.

\bibitem{kim2017data}
S.H. Kim, W.~Whitt, and W.C. Cha.
\newblock A data-driven model of an appointment-generated arrival process at an
  outpatient clinic.
\newblock {\em History}, 2017.

\bibitem{kim2008new}
Y.B. Kim and J.~Park.
\newblock New approaches for inference of unobservable queues.
\newblock In {\em Proceedings of the 40th Conference on Winter Simulation},
  pages 2820--2825. Winter Simulation Conference, 2008.

\bibitem{kim2006congestion}
Y.G. Kim, A.~Shiravi, and P.S. Min.
\newblock Congestion prediction of self-similar network through parameter
  estimation.
\newblock In {\em 2006 IEEE/IFIP Network Operations and Management Symposium
  NOMS 2006}, pages 1--4. IEEE, 2006.

\bibitem{1965Kovalenko}
I.~N Kovalenko.
\newblock In recovering the characteristics of a system from observations of
  the outgoing flow (in russian).
\newblock {\em Dokl. Akad. Nauk SSSR}, 164(5):979--981, 1965.

\bibitem{kraft2009estimating}
S.~Kraft, S.~Pacheco-Sanchez, G.~Casale, and S.~Dawson.
\newblock Estimating service resource consumption from response time
  measurements.
\newblock In {\em Proceedings of the Fourth International ICST Conference on
  Performance Evaluation Methodologies and Tools}, page~48. ICST (Institute for
  Computer Sciences, Social-Informatics and Telecommunications Engineering),
  2009.

\bibitem{krishnasamy2016regret}
S.~Krishnasamy, R.~Sen, R.~Johari, and S.~Shakkottai.
\newblock Regret of queueing bandits.
\newblock {\em Advances in Neural Information Processing Systems}, 29, 2016.

\bibitem{krishnasamy2021learning}
S.~Krishnasamy, R.~Sen, R.~Johari, and S.~Shakkottai.
\newblock Learning unknown service rates in queues: A multiarmed bandit
  approach.
\newblock {\em Operations research}, 69(1):315--330, 2021.

\bibitem{kumar1992average}
A.~Kumar.
\newblock On the average idle time and average queue length estimates in an
  {M}/{M}/1 queue.
\newblock {\em Operations Research Letters}, 12(3):153--157, 1992.

\bibitem{larson1990queue}
R.C. Larson.
\newblock The queue inference engine: deducing queue statistics from
  transactional data.
\newblock {\em Management Science}, pages 586--601, 1990.

\bibitem{larson1991queue}
R.C. Larson.
\newblock The queue inference engine: addendum.
\newblock {\em Management science}, 37(8):1062, 1991.

\bibitem{larson2013queue}
R.C. Larson.
\newblock Queue inference engine.
\newblock In {\em Encyclopedia of Operations Research and Management Science},
  pages 1228--1234. Springer, 2013.

\bibitem{li2019parameter}
C.~Li, H.~Okamura, and T.~Dohi.
\newblock Parameter estimation of {$ M_t$/M/1/K } queueing systems with
  utilization data.
\newblock {\em IEEE Access}, 7:42664--42671, 2019.

\bibitem{li2022hierarchical}
C.~Li, J.~Zheng, H.~Okamura, and T.~Dohi.
\newblock Hierarchical bayesian parameter estimation of queueing systems using
  utilization data.
\newblock {\em International Journal of Performability Engineering}, 18(5),
  2022.

\bibitem{li2022parameter}
C.~Li, J.~Zheng, H.~Okamura, and T.~Dohi.
\newblock Parameter estimation of markovian arrivals with utilization data.
\newblock {\em IEICE Transactions on Communications}, 105(1):1--10, 2022.

\bibitem{doi:10.3141/2623-06}
F.~Li, K.~Tang, J.~Yao, and K.~Li.
\newblock Real-time queue length estimation for signalized intersections using
  vehicle trajectory data.
\newblock {\em Transportation Research Record}, 2623(1):49--59, 2017.

\bibitem{li2013freeway}
L.~Li, X.~Chen, Z.~Li, and L.~Zhang.
\newblock Freeway travel-time estimation based on temporal--spatial queueing
  model.
\newblock {\em Intelligent Transportation Systems, IEEE Transactions on},
  14(3):1536--1541, 2013.

\bibitem{lilliefors1966some}
H.W. Lilliefors.
\newblock Some confidence intervals for queues.
\newblock {\em Operations Research}, pages 723--727, 1966.

\bibitem{lin2021queuing}
Y.~Lin, T.~He, and G.~Pang.
\newblock Queuing network topology inference using passive measurements.
\newblock In {\em 2021 IFIP Networking Conference (IFIP Networking)}, pages
  1--9. IEEE, 2021.

\bibitem{liu2015improving}
H.~Liu, X.~Wu, and P.~Michalopoulos.
\newblock Improving queue size estimation for minnesota's stratified zone
  metering strategy.
\newblock {\em Transportation Research Record: Journal of the Transportation
  Research Board}, 2015.

\bibitem{liu2009real}
H.X. Liu, X.~Wu, W.~Ma, and H.~Hu.
\newblock Real-time queue length estimation for congested signalized
  intersections.
\newblock {\em Transportation research part C: emerging technologies},
  17(4):412--427, 2009.

\bibitem{liu2006adaptive}
X.~Liu, J.~Heo, L.~Sha, and X.~Zhu.
\newblock Adaptive control of multi-tiered web applications using queueing
  predictor.
\newblock In {\em Network Operations and Management Symposium, 2006. NOMS 2006.
  10th IEEE/IFIP}, pages 106--114. IEEE, 2006.

\bibitem{liu2006parameter}
Z.~Liu, L.~Wynter, C.H. Xia, and F.~Zhang.
\newblock Parameter inference of queueing models for it systems using
  end-to-end measurements.
\newblock {\em Performance Evaluation}, 63(1):36--60, 2006.

\bibitem{luo2023queue}
H.~Luo, M.~Deng, J.~Chen, et~al.
\newblock Queue length estimation based on probe vehicle data at signalized
  intersections.
\newblock {\em Journal of Advanced Transportation}, 2023, 2023.

\bibitem{machihara1984carried}
F.~Machihara.
\newblock Carried traffic estimate errors for delay systems.
\newblock {\em Electronics and Communications in Japan (Part I:
  Communications)}, 67(12):49--58, 1984.

\bibitem{mandelbaum1998estimating}
A.~Mandelbaum and S.~Zeltyn.
\newblock Estimating characteristics of queueing networks using transactional
  data.
\newblock {\em Queueing Systems}, 29(1):75--127, 1998.

\bibitem{mandjes2021hypothesis}
M.~Mandjes and L.~Ravner.
\newblock Hypothesis testing for a l{\'e}vy-driven storage system by poisson
  sampling.
\newblock {\em Stochastic Processes and their Applications}, 133:41--73, 2021.

\bibitem{mandjes2005inferring}
M.~Mandjes and R.~van~de Meent.
\newblock Inferring traffic burstiness by sampling the buffer occupancy.
\newblock {\em NETWORKING 2005. Networking Technologies, Services, and
  Protocols; Performance of Computer and Communication Networks; Mobile and
  Wireless Communications Systems}, pages 233--240, 2005.

\bibitem{mandjes2009resource}
M.~Mandjes and R.~van De~Meent.
\newblock Resource dimensioning through buffer sampling.
\newblock {\em IEEE/ACM Transactions on Networking (TON)}, 17(5):1631--1644,
  2009.

\bibitem{mandjes2009queueing}
M.~Mandjes and P.~{\.Z}uraniewski.
\newblock A queueing-based approach to overload detection.
\newblock {\em Network Control and Optimization}, pages 91--106, 2009.

\bibitem{mandjes2011m}
M.~Mandjes and P.~Zuraniewski.
\newblock M/g/[infinity] transience, and its applications to overload
  detection.
\newblock {\em Performance Evaluation}, 2011.

\bibitem{manjunath1996passive}
D.~Manjunath and M.L. Molle.
\newblock Passive estimation algorithms for queueing delays in lans and other
  polling systems.
\newblock In {\em INFOCOM'96. Fifteenth Annual Joint Conference of the IEEE
  Computer Societies. Networking the Next Generation. Proceedings IEEE},
  volume~1, pages 240--247. IEEE, 1996.

\bibitem{manoharan2011markovian}
M.~Manoharan and J.K. Jose.
\newblock {M}arkovian queueing system with random balking.
\newblock {\em OPSEARCH}, 38(3):1--11, 2011.

\bibitem{massey1996estimating}
W.A. Massey, G.A. Parker, and W.~Whitt.
\newblock Estimating the parameters of a nonhomogeneous poisson process with
  linear rate.
\newblock {\em Telecommunication Systems}, 5(2):361--388, 1996.

\bibitem{masuda1995exploiting}
Y.~Masuda.
\newblock Exploiting partial information in queueing systems.
\newblock {\em Operations Research}, 43(3):530--536, 1995.

\bibitem{mccabe2011efficient}
B.P.M. McCabe, G.M. Martin, and D.~Harris.
\newblock Efficient probabilistic forecasts for counts.
\newblock {\em Journal of the Royal Statistical Society: Series B (Statistical
  Methodology)}, 2011.

\bibitem{mcgrath1987subjective}
M.F. Mcgrath, D.~Gross, and N.D. Singpurwalla.
\newblock A subjective bayesian approach to the theory of queues modeling.
\newblock {\em Queueing Systems}, 1(4):317--333, 1987.

\bibitem{mcgrath1987subjectiveII}
M.F. McGrath and N.D. Singpurwalla.
\newblock A subjective bayesian approach to the theory of queues inference
  and information in {M}/{M}/1 queues.
\newblock {\em Queueing Systems}, 1(4):335--353, 1987.

\bibitem{mcvinishconstructing}
R.~McVinish and P.K. Pollett.
\newblock Constructing estimating equations from queue length data.
\newblock {\em preprint}.

\bibitem{mei2019bayesian}
Y.~Mei, W.~Gu, E.C.S. Chung, F.~Li, and K.~Tang.
\newblock A bayesian approach for estimating vehicle queue lengths at
  signalized intersections using probe vehicle data.
\newblock {\em Transportation Research Part C: Emerging Technologies},
  109:233--249, 2019.

\bibitem{meyn2012markov}
S.P. Meyn and R.L. Tweedie.
\newblock {\em Markov chains and stochastic stability}.
\newblock Springer Science \& Business Media, 2012.

\bibitem{mohajerzadeh2015efficient}
A.H. Mohajerzadeh, M.H. Yaghmaee, and A.~Zahmatkesh.
\newblock Efficient data collecting and target parameter estimation in wireless
  sensor networks.
\newblock {\em Journal of Network and Computer Applications}, 57:142--155,
  2015.

\bibitem{mohammadi2012bayesian}
A.~Mohammadi and M.R. Salehi-Rad.
\newblock Bayesian inference and prediction in an {M}/{G}/1 with optional
  second service.
\newblock {\em Communications in Statistics-Simulation and Computation},
  41(3):419--435, 2012.

\bibitem{morales2007bayesian}
J.~Morales, M.~Eugenia~Castellanos, A.M. Mayoral, R.~Fried, and C.~Armero.
\newblock Bayesian design in queues: An application to aeronautic maintenance.
\newblock {\em Journal of statistical planning and inference},
  137(10):3058--3067, 2007.

\bibitem{morozov2016regeneration}
E.~Morozov, R.~Nekrasova, I.~Peshkova, and A.~Rumyantsev.
\newblock A regeneration-based estimation of high performance multiserver
  systems.
\newblock In {\em International Conference on Computer Networks}, pages
  271--282. Springer, 2016.

\bibitem{muddapur1972bayesian}
M.V. Muddapur.
\newblock Bayesian estimates of parameters in some queueing models.
\newblock {\em Annals of the Institute of Statistical Mathematics},
  24(1):327--331, 1972.

\bibitem{muthu1995estimation}
C.~Muthu and V.S. Sampathkumar.
\newblock Estimation of parameters in a particular finite capacity priority
  queueing model.
\newblock {\em Optimization}, 34(4):359--363, 1995.

\bibitem{nam2009estimation}
S.Y. Nam, S.~Kim, and D.K. Sung.
\newblock Estimation of available bandwidth for an {M}/{G}/1 queueing system.
\newblock {\em Applied Mathematical Modelling}, 33(8):3299--3308, 2009.

\bibitem{neal1973theory}
S.R. Neal and A.~Kuczura.
\newblock A theory of traffic measurement errors for loss systems with renewal
  input.
\newblock {\em BSTJ}, 52:967--990, 1973.

\bibitem{nelgabats2012}
Nafna Nelgabats, Yuval Nov, and Gideon Weiss.
\newblock {Sojourn Time Estimation in an M/G/$\infty$ Queue with Partial
  Information}.
\newblock {\em preprint}.

\bibitem{nelson2012stochastic}
B.L. Nelson.
\newblock {\em Stochastic modeling: analysis and simulation}.
\newblock Courier Corporation, 2012.

\bibitem{neuts2005reflections}
M.~Neuts.
\newblock Reflections on statistical methods or complex stochastic systems.
\newblock {\em Modeling Uncertainty}, pages 751--760, 2005.

\bibitem{newell1965approximation}
G.F. Newell.
\newblock Approximation methods for queues with application to the fixed-cycle
  traffic light.
\newblock {\em Siam Review}, 7(2):223--240, 1965.

\bibitem{novak2009determining}
A.~Novak and R.~Watson.
\newblock Determining an adequate probe separation for estimating the arrival
  rate in an {M}/{D}/1 queue using single-packet probing.
\newblock {\em Queueing Systems}, 61(4):255--272, 2009.

\bibitem{nozari1988estimating}
A.~Nozari and W.~Whitt.
\newblock Estimating average production intervals using inventory measurements:
  Little's law for partially observable processes.
\newblock {\em Operations research}, 36(2):308--323, 1988.

\bibitem{ozawa2019stability}
T.~Ozawa.
\newblock Stability condition of a two-dimensional qbd process and its
  application to estimation of efficiency for two-queue models.
\newblock {\em Performance Evaluation}, 130:101--118, 2019.

\bibitem{pakes1971serial}
A.G. Pakes.
\newblock The serial correlation coefficients of waiting times in the
  stationary {GI/M/1} queue.
\newblock {\em The Annals of Mathematical Statistics}, 42(5):1727--1734, 1971.

\bibitem{park2007choice}
J.~Park.
\newblock On the choice of an auxiliary function in the {$M/G/\infty$}
  estimation.
\newblock {\em Computational Statistics \& Data Analysis}, 51(12):5477--5482,
  2007.

\bibitem{park2011analysis}
J.~Park, Y.B. Kim, and T.R. Willemain.
\newblock Analysis of an unobservable queue using arrival and departure times.
\newblock {\em Computers \& Industrial Engineering}, 61(3):842--847, 2011.

\bibitem{paschalidis2001estimation}
I.C. Paschalidis and S.~Vassilaras.
\newblock On the estimation of buffer overflow probabilities from measurements.
\newblock {\em Information Theory, IEEE Transactions on}, 47(1):178--191, 2001.

\bibitem{pichitlamken2003modelling}
Pichitlamken, Deslauriers, and Avramidis.
\newblock Modelling and simulation of a telephone call center.
\newblock In {\em Proceedings of the 2003 Winter Simulation Conference, 2003.},
  volume~2, pages 1805--1812. IEEE, 2003.

\bibitem{pickands1997estimation}
J.~Pickands~III and R.A. Stine.
\newblock Estimation for an {M/G/$\infty$} queue with incomplete information.
\newblock {\em Biometrika}, 84(2):295--308, 1997.

\bibitem{pin2010statistical}
F.~Pin, D.~Veitch, and B.~Kauffmann.
\newblock Statistical estimation of delays in a multicast tree using
  accelerated em.
\newblock {\em Queueing Systems}, 66:369--412, 2010.

\bibitem{pitts1994nonparametric}
S.M. Pitts.
\newblock Nonparametric estimation of the stationary waiting time distribution
  function for the gi/g/1 queue.
\newblock {\em The Annals of Statistics}, pages 1428--1446, 1994.

\bibitem{7850972}
N.~Polson and V.~Sokolov.
\newblock Bayesian particle tracking of traffic flows.
\newblock {\em IEEE Transactions on Intelligent Transportation Systems},
  19(2):345--356, 2018.

\bibitem{prabhu2012stochastic}
N.U. Prabhu.
\newblock {\em Stochastic storage processes: queues, insurance risk, dams, and
  data communication}, volume~15.
\newblock Springer Science \& Business Media, 2012.

\bibitem{prieger2005estimation}
J.E. Prieger.
\newblock Estimation of a simple queuing system with units-in-service and
  complete data.
\newblock {\em Working Papers}, 2005.

\bibitem{quinino2016bayesian}
R.C. Quinino and F.R.B. Cruz.
\newblock Bayesian sample sizes in an {M/M/1} queueing systems.
\newblock {\em The International Journal of Advanced Manufacturing Technology},
  pages 1--8, 2016.

\bibitem{quinino2017bayesian}
R.C. Quinino and F.R.B. Cruz.
\newblock Bayesian sample sizes in an {M/M/1} queueing systems.
\newblock {\em The International Journal of Advanced Manufacturing Technology},
  88(1):995--1002, 2017.

\bibitem{ramalhoto1987some}
M.F. Ramalhoto.
\newblock Some statistical problems in random translations of stochastic point
  processes.
\newblock {\em Annals of Operations Research}, 8(1):229--242, 1987.

\bibitem{ramirez2008bayesian}
P.~Ramirez, R.E. Lillo, and M.P. Wiper.
\newblock Bayesian analysis of a queueing system with a long-tailed arrival
  process.
\newblock {\em Communications in Statistics--Simulation and Computation},
  37(4):697--712, 2008.

\bibitem{ramirez2008bayesian-thesis}
J.~Ram{\'\i}rez~Cobo.
\newblock Bayesian modelling of stochastic processes in teletraffic and
  finance.
\newblock 2008.

\bibitem{ramirez2010bayesian}
P.~Ramirez-Cobo, R.E. Lillo, S.~Wilson, and M.P. Wiper.
\newblock Bayesian inference for double pareto lognormal queues.
\newblock {\em The Annals of Applied Statistics}, 4(3):1533--1557, 2010.

\bibitem{ravner2022queue}
L.~Ravner.
\newblock Queue input estimation from discrete workload observations.
\newblock {\em Queueing Systems}, 100(3-4):541--543, 2022.

\bibitem{2019Ravner}
L.~Ravner, O.~Boxma, and M.~Mandjes.
\newblock {Estimating the input of a Lévy-driven queue by Poisson sampling of
  the workload process}.
\newblock {\em Bernoulli}, 25(4B):3734 -- 3761, 2019.

\bibitem{ravner2023estimating}
L.~Ravner and J.~Wang.
\newblock Estimating customer delay and tardiness sensitivity from periodic
  queue length observations.
\newblock {\em Queueing Systems}, 103(3-4):241--274, 2023.

\bibitem{ren2012bayes}
H.~Ren and G.~Wang.
\newblock Bayes estimation of traffic intensity in {M/M/1} queue under a
  precautionary loss function.
\newblock {\em Procedia Engineering}, 29:3646--3650, 2012.

\bibitem{ren2012bayesB}
H.P. Ren and J.P. Li.
\newblock Bayes estimation of traffic intensity in {M/M/1} queue under a new
  weighted square error loss function.
\newblock {\em Advanced Materials Research}, 485:490--493, 2012.

\bibitem{reynolds1973estimating}
J.F. Reynolds.
\newblock On estimating the parameters of a birth-death process.
\newblock {\em Australian Journal of Statistics}, 15(1):35--43, 1973.

\bibitem{rodrigo2006estimators}
A.~Rodrigo.
\newblock Estimators of the retrial rate in m/g/1 retrial queues.
\newblock {\em ASIA PACIFIC JOURNAL OF OPERATIONAL RESEARCH}, 23(2):193, 2006.

\bibitem{rodrigo1999large}
A.~Rodrigo and M.~Vazquez.
\newblock Large sample inference in retrial queues.
\newblock {\em Mathematical and Computer Modelling}, 30(3-4):197--206, 1999.

\bibitem{rodrigues1998note}
J.~Rodrigues and J.G. Leite.
\newblock A note on bayesian analysis in {M}/{M}/1 queues derived from
  confidence intervals.
\newblock {\em Statistics: A Journal of Theoretical and Applied Statistics},
  31(1):35--42, 1998.

\bibitem{ross2005estimating}
A.M. Ross and J.G. Shanthikumar.
\newblock Estimating effective capacity in erlang loss systems under
  competition.
\newblock {\em Queueing Systems}, 49(1):23--47, 2005.

\bibitem{ross2007estimation}
J.V. Ross, T.~Taimre, and P.K. Pollett.
\newblock Estimation for queues from queue length data.
\newblock {\em Queueing Systems}, 55(2):131--138, 2007.

\bibitem{ross1970identifiability}
S.M. Ross.
\newblock Identifiability in gi/g/k queues.
\newblock {\em Journal of Applied Probability}, 7(3):776--780, 1970.

\bibitem{rubin1990single}
G.~Rubin and D.S. Robson.
\newblock A single server queue with random arrivals and balking: confidence
  interval estimation.
\newblock {\em Queueing Systems}, 7(3):283--306, 1990.

\bibitem{schruben1982some}
L.~Schruben and R.~Kulkarni.
\newblock Some consequences of estimating parameters for the {M}/{M}/1 queue.
\newblock {\em Operations Research Letters}, 1(2):75--78, 1982.

\bibitem{schweer2015nonparametric}
S.~Schweer and C.~Wichelhaus.
\newblock Nonparametric estimation of the service time distribution in the
  discrete-time {$GI/G/\infty$} queue with partial information.
\newblock {\em Stochastic Processes and their Applications}, 125(1):233--253,
  2015.

\bibitem{senderovich2014queue}
A.~Senderovich, M.~Weidlich, A.~Gal, and A.~Mandelbaum.
\newblock Queue mining--predicting delays in service processes.
\newblock In {\em International Conference on Advanced Information Systems
  Engineering}, pages 42--57. Springer, 2014.

\bibitem{senderovich2015discovering}
Arik Senderovich, S~Leemans, Shahar Harel, Avigdor Gal, Avishai Mandelbaum, and
  W~van~der Aalst.
\newblock Discovering queues from event logs with varying levels of
  information.
\newblock {\em Lect Notes Bus Inf (forthcoming)}, 2015.

\bibitem{senderovich2015queue}
Arik Senderovich, Matthias Weidlich, Avigdor Gal, and Avishai Mandelbaum.
\newblock Queue mining for delay prediction in multi-class service processes.
\newblock {\em Information Systems}, 2015.

\bibitem{sengupta1989markov}
B.~Sengupta.
\newblock Markov processes whose steady state distribution is
  matrix-exponential with an application to the {GI/PH/1} queue.
\newblock {\em Advances in Applied Probability}, pages 159--180, 1989.

\bibitem{sharma1999estimatingB}
V.~Sharma.
\newblock Estimating traffic intensities at different nodes in networks via a
  probing stream.
\newblock In {\em Global Telecommunications Conference, 1999. GLOBECOM'99},
  volume~1, pages 374--380. IEEE, 1999.

\bibitem{sharma1998estimating}
V.~Sharma and R.~Mazumdar.
\newblock Estimating traffic parameters in queueing systems with local
  information.
\newblock {\em Performance Evaluation}, 32(3):217--230, 1998.

\bibitem{singh2021bayesian}
S.K. Singh, S.K. Acharya, F.R.B. Cruz, and R.C. Quinino.
\newblock Bayesian sample size determination in a single-server deterministic
  queueing system.
\newblock {\em Mathematics and Computers in Simulation}, 187:17--29, 2021.

\bibitem{singh2022bayesian}
S.K. Singh, S.K. Acharya, F.R.B. Cruz, and R.C. Quinino.
\newblock Bayesian inference and prediction in an queueing system.
\newblock {\em Communications in Statistics-Theory and Methods}, pages 1--21,
  2022.

\bibitem{singpurwalla1992discussion}
N.D. Singpurwalla.
\newblock Discussion of Thiruvaiyaru and Basawa's empirical bayes estimation
  for queueing systems and networks.
\newblock {\em Queueing Systems}, 11(3):203--206, 1992.

\bibitem{smith1953distribution}
W.L. Smith.
\newblock On the distribution of queueing times.
\newblock In {\em Mathematical Proceedings of the Cambridge Philosophical
  Society}, volume~49, pages 449--461. Cambridge Univ Press, 1953.

\bibitem{sohn1996empirical}
S.Y. Sohn.
\newblock Empirical bayesian analysis for traffic intensity: {M}/{M}/1 queues
  with covariates.
\newblock {\em Queueing Systems}, 22(3):383--401, 1996.

\bibitem{sohn1996influence}
S.Y. Sohn.
\newblock Influence of a prior distribution on traffic intensity estimation
  with covariates.
\newblock {\em Journal of Statistical Computation and Simulation},
  55(3):169--180, 1996.

\bibitem{sohn2002robust}
S.Y. Sohn.
\newblock Robust design of server capability in {M}/{M}/1 queues with both
  partly random arrival and service rates.
\newblock {\em Computers \& Operations Research}, 29(5):433--440, 2002.

\bibitem{sousa2011suitability}
M.E. Sousa-Vieira.
\newblock Suitability of the {M/G/$\infty$} process for modeling scalable h.
  264 video traffic.
\newblock {\em Analytical and Stochastic Modeling Techniques and Applications},
  pages 149--158, 2011.

\bibitem{srinivas2015ml}
V.~Srinivas and B.K. Kale.
\newblock {ML} and {UMVU} estimation in the {M/D/1 }queuing system.
\newblock {\em Communications in Statistics-Theory and Methods},
  (just-accepted), 2015.

\bibitem{srinivas2011estimation}
V.~Srinivas, S.S. Rao, and B.K. Kale.
\newblock Estimation of measures in {{M}/{M}/1} queue.
\newblock {\em Communications in Statistics-Theory and Methods},
  40(18):3327--3336, 2011.

\bibitem{subba1986large}
S.~Subba~Rao and K.~Harishchandra.
\newblock On a large sample test for the traffic intensity in {GI}/{G}/s queue.
\newblock {\em Naval Research Logistics Quarterly}, 33(3):545--550, 1986.

\bibitem{sutarto2015modeling}
H.Y. Sutarto and E.~Joelianto.
\newblock Modeling, identification, estimation, and simulation of urban traffic
  flow in jakarta and bandung.
\newblock {\em Journal of Mechatronics, Electrical Power, and Vehicular
  Technology}, 6(1):57--66, 2015.

\bibitem{sutarto2017developing}
H.Y. Sutarto, E.~Joelianto, and T.A. Nugroho.
\newblock Developing a stochastic model of queue length at a signalized
  intersection.
\newblock {\em International Journal on Advanced Science, Engineering and
  Information Technology}, 7:2183--2188, 2017.

\bibitem{sutton2008probabilistic}
C.~Sutton and M.I. Jordan.
\newblock Probabilistic inference in queueing networks.
\newblock In {\em Proceedings of the Third conference on Tackling computer
  systems problems with machine learning techniques}, pages 6--6. USENIX
  Association, 2008.

\bibitem{sutton2010learning}
C.~Sutton and M.I. Jordan.
\newblock Learning and inference in queueing networks.
\newblock 2010.

\bibitem{sutton2011bayesian}
C.~Sutton and M.I. Jordan.
\newblock Bayesian inference for queueing networks and modeling of internet
  services.
\newblock {\em The Annals of Applied Statistics}, 5(1):254--282, 2011.

\bibitem{suyama2018simple}
E.~Suyama, R.C. Quinino, and F.R.B. Cruz.
\newblock Simple and yet efficient estimators for markovian multiserver queues.
\newblock {\em Mathematical Problems in Engineering}, 2018, 2018.

\bibitem{tan2020fuzing}
C.~Tan, L.~Liu, H.~Wu, Y.~Cao, and K.~Tang.
\newblock Fuzing license plate recognition data and vehicle trajectory data for
  lane-based queue length estimation at signalized intersections.
\newblock {\em Journal of Intelligent Transportation Systems}, 24(5):449--466,
  2020.

\bibitem{tan2019cycle}
C.~Tan, J.~Yao, K.~Tang, and J.~Sun.
\newblock Cycle-based queue length estimation for signalized intersections
  using sparse vehicle trajectory data.
\newblock {\em IEEE Transactions on Intelligent Transportation Systems},
  22(1):91--106, 2019.

\bibitem{thiagarajan1979statistical}
T.R. Thiagarajan and C.M. Harris.
\newblock {Statistical tests for exponential service from M/G/1 waiting-time
  data}.
\newblock {\em Naval Research Logistics Quarterly}, 26(3):511--520, 1979.

\bibitem{thiruvaiyaru1992empirical}
D.~Thiruvaiyaru and I.V. Basawa.
\newblock Empirical bayes estimation for queueing systems and networks.
\newblock {\em Queueing Systems}, 11(3):179--202, 1992.

\bibitem{thiruvaiyaru1991estimation}
D.~Thiruvaiyaru, I.V. Basawa, and U.N. Bhat.
\newblock Estimation for a class of simple queueing networks.
\newblock {\em Queueing Systems}, 9(3):301--312, 1991.

\bibitem{toyoizumi1997sengupta}
H.~Toyoizumi.
\newblock Sengupta's invariant relationship and its application to waiting time
  inference.
\newblock {\em Journal of Applied Probability}, 34(3):795--799, 1997.

\bibitem{van2020estimation}
C.N. Van~Phu and N.~Farhi.
\newblock Estimation of urban traffic state with probe vehicles.
\newblock {\em IEEE Transactions on Intelligent Transportation Systems},
  22(5):2797--2808, 2020.

\bibitem{veeger2011predicting}
C.P.L. Veeger, L.F.P. Etman, E.~Lefeber, I.J.B.F. Adan, J.~Van~Herk, and
  JE~Rooda.
\newblock Predicting cycle time distributions for integrated processing
  workstations: an aggregate modeling approach.
\newblock {\em IEEE Transactions on Semiconductor Manufacturing},
  24(2):223--236, 2011.

\bibitem{vorobeychikov2015cusum}
S.~Vorobeychikov.
\newblock Cusum algorithms for parameter estimation in queueing systems with
  jump intensity of the arrival process.
\newblock In {\em Information Technologies and Mathematical Modelling-Queueing
  Theory and Applications: 14th International Scientific Conference, ITMM 2015,
  named after AF Terpugov, Anzhero-Sudzhensk, Russia, November 18-22, 2015,
  Proceedings}, volume 564, page 275. Springer, 2015.

\bibitem{walrand1981filtering}
J.~Walrand.
\newblock Filtering formulas and the ./{M}/1 queue in a quasireversible
  network.
\newblock {\em Stochastics: An International Journal of Probability and
  Stochastic Processes}, 6(1):1--22, 1981.

\bibitem{walrand1988introduction}
J.~Walrand.
\newblock {\em An introduction to queueing networks}, volume~21.
\newblock Prentice Hall Englewood Cliffs, NJ, 1988.

\bibitem{walton2021learning}
Neil Walton and Kuang Xu.
\newblock Learning and information in stochastic networks and queues.
\newblock In {\em Tutorials in Operations Research: Emerging Optimization
  Methods and Modeling Techniques with Applications}, pages 161--198. INFORMS,
  2021.

\bibitem{wang2021calibrating}
R.~Wang and H.~Honnappa.
\newblock Calibrating infinite server queueing models driven by cox processes.
\newblock In {\em 2021 Winter Simulation Conference (WSC)}, pages 1--12. IEEE,
  2021.

\bibitem{wang2020combining}
S.~Wang, W.~Huang, and H.K. Lo.
\newblock Combining shockwave analysis and bayesian network for traffic
  parameter estimation at signalized intersections considering queue spillback.
\newblock {\em Transportation Research Part C: Emerging Technologies},
  120:102807, 2020.

\bibitem{wang2006maximum}
T.Y. Wang, J.C. Ke, K.H. Wang, and S.C. Ho.
\newblock Maximum likelihood estimates and confidence intervals of an
  {M}/{M}/{R} queue with heterogeneous servers.
\newblock {\em Mathematical Methods of Operations Research}, 63(2):371--384,
  2006.

\bibitem{wang2015maximum}
W.~Wang and G.~Casale.
\newblock Maximum likelihood estimation of closed queueing network demands from
  queue length data.
\newblock {\em ACM SIGMETRICS Performance Evaluation Review}, 43(2):45--47,
  2015.

\bibitem{wang2015filling}
W.~Wang, J.F. P{\'e}rez, and G.~Casale.
\newblock Filling the gap: a tool to automate parameter estimation for software
  performance models.
\newblock In {\em Proceedings of the 1st International Workshop on
  Quality-Aware DevOps}, pages 31--32. ACM, 2015.

\bibitem{warfield1984application}
R.E. Warfield and G.A. Foers.
\newblock Application of bayesian methods to teletraffic measurement and
  dimensioning.
\newblock {\em Australian Telecommunications Research}, 18:51--58, 1984.

\bibitem{warfield1985application}
R.E. Warfield and G.A. Foers.
\newblock Application of bayesian teletraffic measurement to systems with
  queueing or repeated attempts.
\newblock In {\em Proceedings of the Eleventh International Teletraffic
  Congress}, 1985.

\bibitem{weerasinghe2013abandonment}
A.~Weerasinghe and A.~Mandelbaum.
\newblock Abandonment versus blocking in many-server queues: asymptotic
  optimality in the {QED} regime.
\newblock {\em Queueing Systems}, 75(2-4):279--337, 2013.

\bibitem{whitt2012fitting}
W.~Whitt.
\newblock Fitting birth-and-death queueing models to data.
\newblock {\em Statistics and Probability Letters}, 82:998--1004, 2012.

\bibitem{whitt2015many}
W.~Whitt.
\newblock Many-server limits for periodic infinite-server queues.
\newblock {\em Columbia University}, 2015.

\bibitem{whitt2017data}
W.~Whitt and X.~Zhang.
\newblock A data-driven model of an emergency department.
\newblock {\em Operations Research for Health Care}, 12:1--15, 2017.

\bibitem{wiler2013emergency}
J.L. Wiler, E.~Bolandifar, R.T. Griffey, R.F. Poirier, and T.~Olsen.
\newblock An emergency department patient flow model based on queueing theory
  principles.
\newblock {\em Academic Emergency Medicine}, 20(9):939--946, 2013.

\bibitem{wiper1998bayesian}
M.P. Wiper.
\newblock Bayesian analysis of {$E_r$}/{M}/1 and {$E_r$}/{M}/c queues.
\newblock {\em Journal of Statistical Planning and Inference}, 69(1):65--79,
  1998.

\bibitem{wolff1965problems}
R.W. Wolff.
\newblock Problems of statistical inference for birth and death queuing models.
\newblock {\em Operations Research}, pages 343--357, 1965.

\bibitem{woodside1984optimal}
C.M. Woodside, D.A. Stanford, and B.~Pagurek.
\newblock Optimal prediction of queue lengths and delays in gi/m/m multiserver
  queues.
\newblock {\em Operations research}, pages 809--817, 1984.

\bibitem{xu2010hypothesis}
X.~Xu, Q.~Zhang, and X.~Ding.
\newblock Hypothesis testing and confidence regions for the mean sojourn time
  of an {M}/{M}/1 queueing system.
\newblock {\em Communications in Statistics Theory and Methods}, 40(1):28--39,
  2010.

\bibitem{yom2014erlang}
G.B. Yom-Tov and A.~Mandelbaum.
\newblock Erlang-r: A time-varying queue with reentrant customers, in support
  of healthcare staffing.
\newblock {\em Manufacturing \& Service Operations Management}, 16(2):283--299,
  2014.

\bibitem{yu2015hidden}
S.Z. Yu.
\newblock {\em Hidden Semi-Markov Models: Theory, Algorithms and Applications}.
\newblock Morgan Kaufmann, 2015.

\bibitem{zammit2016joint}
L.C. Zammit, S.G. Fabri, and K.~Scerri.
\newblock Joint state and parameter estimation for a macro traffic junction
  model.
\newblock In {\em Control and Automation (MED), 2016 24th Mediterranean
  Conference on}, pages 1152--1157. IEEE, 2016.

\bibitem{8317685}
L.C. Zammit, S.G. Fabri, and K.~Scerri.
\newblock Online state and multidimensional parameter estimation for a
  macroscopic model of a traffic junction.
\newblock In {\em 2017 IEEE 20th International Conference on Intelligent
  Transportation Systems (ITSC)}, pages 1--6, 2017.

\bibitem{zhan2015lane}
X.~Zhan, R.~Li, and S.V. Ukkusuri.
\newblock Lane-based real-time queue length estimation using license plate
  recognition data.
\newblock {\em Transportation Research Part C: Emerging Technologies},
  57:85--102, 2015.

\bibitem{8759936}
H.~Zhang, H.X. Liu, P.~Chen, G.~Yu, and Y.~Wang.
\newblock Cycle-based end of queue estimation at signalized intersections using
  low-penetration-rate vehicle trajectories.
\newblock {\em IEEE Transactions on Intelligent Transportation Systems},
  21(8):3257--3272, 2020.

\bibitem{zhang2002workload}
L.~Zhang, C.H. Xia, M.S. Squillante, and W.N. Mills~III.
\newblock Workload service requirements analysis: A queueing network
  optimization approach.
\newblock In {\em Modeling, Analysis and Simulation of Computer and
  Telecommunications Systems, 2002. MASCOTS 2002. Proceedings. 10th IEEE
  International Symposium on}, pages 23--32. IEEE, 2002.

\bibitem{zhang2010confidence}
Q.~Zhang and X.~Xu.
\newblock Confidence intervals of performance measures for an {M/G/1} queueing
  system.
\newblock {\em Communications in Statistics Simulation and
  Computation{\textregistered}}, 39(3):501--516, 2010.

\bibitem{zhang2016generalized}
Q.~Zhang, X.~Xu, and S.~Mi.
\newblock A generalized p-value approach to inference on the performance
  measures of an {$M/E_k/1$} queueing system.
\newblock {\em Communications in Statistics-Theory and Methods},
  45(8):2256--2267, 2016.

\bibitem{zhao2020traffic}
Y.~Zhao.
\newblock {\em Traffic State Estimation Using Probe Vehicle Data}.
\newblock PhD thesis, 2020.

\bibitem{zhao2021hidden}
Y.~Zhao, S.~Shen, and H.X. Liu.
\newblock A hidden markov model for the estimation of correlated queues in
  probe vehicle environments.
\newblock {\em Transportation Research Part C: Emerging Technologies},
  128:103128, 2021.

\bibitem{zhao2021maximum}
Y.~Zhao, W.~Wong, J.~Zheng, and H.X. Liu.
\newblock Maximum likelihood estimation of probe vehicle penetration rates and
  queue length distributions from probe vehicle data.
\newblock {\em IEEE Transactions on Intelligent Transportation Systems},
  23(7):7628--7636, 2021.

\bibitem{zhao2019various}
Y.~Zhao, J.~Zheng, W.~Wong, X.~Wang, Y.~Meng, and H.X. Liu.
\newblock Various methods for queue length and traffic volume estimation using
  probe vehicle trajectories.
\newblock {\em Transportation Research Part C: Emerging Technologies},
  107:70--91, 2019.

\bibitem{zheng2000some}
S.~Zheng and A.F. Seila.
\newblock Some well-behaved estimators for the {M/M/1} queue.
\newblock {\em Operations Research Letters}, 26(5):231--235, 2000.

\bibitem{zhong2022learning}
Y.~Zhong, J.R. Birge, and A.~Ward.
\newblock Learning the scheduling policy in time-varying multiclass many server
  queues with abandonment.
\newblock {\em Available at SSRN}, 2022.

\bibitem{zuraniewski2010empirical}
P.~Zuraniewski, M.~Mandjes, and M.~Mellia.
\newblock Empirical assessment of voip overload detection tests.
\newblock In {\em Next Generation Internet (NGI), 2010 6th EURO-NF Conference
  on}, pages 1--8. IEEE, 2010.

\end{thebibliography}
				
			\end{document}